\documentclass[a4paper]{article}
\usepackage{a4wide}
\usepackage[utf8]{inputenc}
\usepackage[english]{babel}
\usepackage{amssymb}
\usepackage{amsfonts}
\usepackage{amsmath}
\usepackage{amsthm}
\usepackage{graphicx}
\usepackage{hyperref}
\usepackage{subfigure}

\theoremstyle{remark}

\theoremstyle{definition}

\newtheorem{theor}{Theorem}
\newtheorem{cor}{Corollary}
\newtheorem{ex}{Example}
\newtheorem{prop}{Proposition}

\newtheorem{lem}{Lemma}

\theoremstyle{remark}
\newtheorem{rem}
{R~e~m~a~r~k}

\begin{document}

\begin{center}
\LARGE{Universal local linear kernel estimators in nonparametric regression}
\end{center}
\begin{center}
\textit{Yuliana Linke $^{1}$, Igor Borisov $^{1}$,
 Pavel Ruzankin $^{1}$, Vladimir Kutsenko $^{2,3}$ , Elena Yarovaya $^{2,3}$ and Svetlana Shalnova $^{3}$}
\end{center}

\begin{center}
$^{1}$  Sobolev Institute of Mathematics, Novosibirsk, Russia;\\

$^{2}$  Department of Probability Theory, Lomonosov Moscow State University, Moscow, Russia;\\

$^{3}$  Department of Epidemiology of Noncommunicable Diseases, National Medical Research Center for Therapy and Preventive Medicine, Moscow, Russia;
\end{center}
\medskip

\begin{abstract}
New local linear estimators are proposed for a wide class of nonparametric regression models.
The estimators are uniformly consistent regardless of satisfying
 traditional conditions of depen\-dence of design elements.
 The estimators are the solutions of a specially weighted least-squares method.
The design can be fixed or random and does not need to meet
 classical regularity or independence conditions.
 As an application, several estimators are constructed for the
 mean of dense functional data.
 The theoretical results of the study are illustrated by simulations.
 An example of processing real medical data from the epidemiological cross-sectional study ESSE-RF  is included.
 We compare the new estimators with the estimators best known for such studies.

\medskip
\textbf{Keywords:} nonparametric regression; kernel estimator; local linear estimator; uniform consistency; fixed design; random design; dependent design elements; mean of dense functional data; epidemiological research.
\medskip

\textbf{2020 Mathematics subject classification:} 62G08
\end{abstract}

\section{Introduction}

In this paper, we consider a nonparametric regression model, where bivariate 
observations \[\{(X_1,z_1),\ldots,(X_n,z_n)\}\] satisfy
the following equations:
 \begin{equation}  \label{1}
 X_i=f({z}_i)+\varepsilon_i,
 \qquad i=1,\ldots,n,
 \end{equation}
where
$f(t)\equiv f(\omega, t)$, $t\in [0,1]$, is an unknown random function (process) which is continuous almost surely,
the design  $\{z_i;\,i=1,\ldots,n\}$  consists of a set of observable random variables with possibly unknown distributions
lying
in~$[0,1]$, the design points are not necessarily independent or identically distributed.
We will consider the design as a triangular array, i.e., the random variables $\{z_i;\; i = 1, \ldots, n\}$ may depend on $n$. In particular, this scheme includes
regression models with fixed design. The random regression function $f(t)$ is not supposed to be design independent. We will give below some fairly standard conditions for
the regression analysis
on the random errors $\{\varepsilon_i;\, i = 1, \ldots, n \}$. In particular, they are supposed to be centered, not necessarily independent
or identically distributed.

The paper is devoted to  constructing uniformly consistent estimators for the regression function $f (t)$ under minimal assumptions
on the correlation of design points.

The most popular kernel estimation procedures in the classical case of nonrandom regression function
are apparently related
with the estimators of Nadaray--Watson, Priest\-ley--Zhao, Gasser--M\"{u}ller, local polynomial estimators, as well as their modifications (e.g., see
 \cite{1996-Fa}, \cite{2003-FY},    \cite{2002-Gy},  \cite{1990-Ha},  \cite{1988-M}).
We are primarily interested in the dependence conditions of design elements $\{{z}_i\}$. In this regard, a huge number of publications in the field
of nonparametric regression can be conditionally divided into two groups.
We will classify papers with a random design to the first one, and to the second one -- with a fixed design.

In the papers dealing with random design, either independent and identically distributed quantities are considered or, as a rule, stationary sequences of observations
that satisfy one or another known form of dependence. In particular, various types of mixing conditions, schemes of moving averages,
associated random variables, Markov or martingale properties, and so on have been used. In this regard, we note, for example, the papers
\cite{2003-CD},  \cite{1979-Dev}, \cite{1990-Ga}, \cite{2002-Gy},    \cite{2008-Ha},    \cite{1984-HL}, \cite{2016-HL}, \cite{2001-JM}, \cite{2011-K1}, \cite{1989-Lie}, \cite{2010-LJ}, \cite{2016-LYS}--\cite{2005-Ma},  \cite{1997-Mu}, \cite{1970-Nad}, \cite{1990-Ro}, \cite{2013-SX}.
In the recent papers \cite{2012-CGL}, \cite{2007-KMT},
\cite{2016-LW},
 and \cite{2014-WC}, nonstationary sequences of design elements with one or another special type of dependence are considered
(Markov chains, autoregression, partial sums of moving averages, etc.).
In the case of fixed design,
  in the overwhelming majority of works, certain conditions for the regularity of the design are assumed
  (e.g., see \cite{2020-BBL}--\cite{2001-BF}, \cite{2007-GRT}, \cite{2008-Ha}, \cite{1984-HL},
  \cite{2018-TXW},
\cite{1994-Wu},   \cite{2020-ZZ}).
So,   the
nonrandom design points $z_i$ are most often given by the formula $z_i = g (i/n) + o (1/n)$ with some function $g$ of bounded variation, where the error $ o (1 / n) $
is uniform in all $i = 1,\ldots, n$. If $g$ is linear then we get a so-called {\it equidistant} design.
Another version of the regularity condition is the relation $\max\nolimits_{i \leq n} (z_i-z_{i-1}) = O(1/n)$ (here it is assumed that the design elements
ranged in increasing order).

The problem of uniform approximation of a regression function has been studied by many authors (e.g., see
\cite{1979-Dev}, \cite{2005-Ei},  \cite{2007-GRT},  \cite{2008-Ha}, \cite{1984-HL}, \cite{1993-I},  \cite{2005-LJ}, \cite{1989-Lie}, \cite{2016-LYS},  \cite{1982-MS}, \cite{1970-Nad}, \cite{2013-SX}, \cite{2014-WC}, \cite{ZLL-2018},
  and the references there).

In connection with studying the random regression function $f(t)$, we note, for example, the papers
 \cite{2006-Ha}, \cite{KR-2017},  \cite{2010-Li}, \cite{LJ-2020},
 \cite{Y-2007}--\cite{ZLL-2018}
 where the mean and covariance functions of the random regression function~$f$ are estimated in the case when, for $N$
 indepen\-dent copies $f_1,\ldots,f_N$ of the function $f$, noisy values of each of these trajectories are observed for some collection of design elements (the design can be either common to all trajectories or different from series to series).
 Estimation of the mean and covariance functions is an actively developing area of nonparametric estimation, especially in the last couple of decades, which is both of independent interest and plays an important role for some subsequent analysis  of the random process $f$
(e.g., see \cite{HE-2015}, \cite{KR-2017},  \cite{2010-Li}, \cite{Mu-2005},  \cite{WCM-2016}, \cite{ZW-2016}). We consider one of the variants of this problem as an application of the main result.

The purpose of this article is to construct estimators that are uniformly consistent (in the sense of convergence in probability) not only
in the above-mentioned
 review of cases of dependence, but also for significantly different correlations of observations when the conditions of ergodicity
or stationarity are not satisfied,
as well as the classical mixing conditions and other well-known dependence restrictions.
Note that the proposed estimators belong to the class of local linear kernel estimators, but with some different weights than in the classical version.
In this case, instead of the original observations, we consider their concomitants associated with the variational series based on the design observations, and their spacings are taken as the additional weights for the corresponding weighted least-square method generating the above-mentioned new estimators.
It is important to emphasize that these estimators have the property of universality regarding the nature of dependence of observations:
the design can be either fixed and
not necessarily regular, or random, while not necessarily satisfying the traditional correlation conditions.
   In particular, the only condition for design points that guarantees the uniform consistency of new estimators is the condition
for dense filling of the domain of definition of the regression function.
In our opinion, this condition is very clear and
in fact, it is necessary to restore the function on the area of defining design elements.
Previously, similar ideas were implemented in \cite{2020} for slightly different evaluations (in detail, see Sec. 4).
Similar conditions for design elements were also used in \cite{2022} and \cite{2021} in nonparametric regression, and in \cite{2019}--\cite{2018-eng}
--- in nonlinear regression.

The paper has the following structure. Section 2 contains the main results,
Section 3 discusses
the problem of estimating the mean function of a stochastic process.
Comparison of the universal local linear estimators with some known ones
is given in Section~4.
Section~5  contains some results of computer simulation. In Section~6, we compare the results of using the new universal local linear estimators
with the most common approaches of data analysis based on the epidemiological research ESSE-RF.
In Section 7, we briefly summarize the results of the study.
The proofs of the results from Sections 2--4 are referred to Section~8.

\section{Main results}

We need a number of assumptions.

$({\bf D})$ {\it The observations $X_1,\ldots, X_n$ are represented in the form {\rm(\ref{1})}, where 
the unknown random regression function
$f:[0,1]\to {\mathbb R}$ is continuous almost surely.
The design points  $\{z_i;\,i=1,\ldots,n\}$
are a set of observable random variables with values in $[0,1]$, having, generally speaking, unknown distributions, not necessarily independent or equally distributed.
Moreover, the random variables $\{z_i; \,i = 1, \ldots,n\} $ may depend on $n$, i.e., can be considered as an array of design observations. The random function $f(t)$ may be design dependent.}

$({\bf E})$ {\it  For all $n \geq 1$, the unobservable random errors $\{\varepsilon_i; \, i = 1, \ldots,n\}$ satisfy with probability $1$
the following conditions for all $i, j \le n$ and $i \neq j$$:$
\begin{equation}\label{w2}
 {\mathbb E}_{{\cal F}_n}\varepsilon_i=0,\qquad
 \sup_{i\leq n}{\mathbb E}_{{\cal F}_n}\varepsilon^2_i\le \sigma^2,
\qquad {\mathbb E}_{{\cal F}_n}\varepsilon_i\varepsilon_j= 0,
\end{equation}
where the constant $\sigma^2>0$ may be unknown and does not depend on $n$, the symbol ${\mathbb E}_{{\cal F}_n}$ stands for the conditional expectation given the $\sigma$-field generated both by the paths of the random process $f(\cdot)$ and by the random variables $\{z_{i}; \,i=1,\ldots,n\}$.
}

$({\bf K})$
{\it  A kernel $K(t)$,  $t\in \mathbb R$, is equal to zero outside the interval $[-1,1]$ and is the density of a symmetric distribution with the support in $[-1,1]$, i.e.,
$K(t)\ge 0 $, $K(t)=K(-t)$ for all  $t\in [-1,1]$, and  $\int\nolimits_{-1}^1 K(t) dt =1$.
We assume that the function  $K(t)$ satisfies the Lipschitz condition with constant  $1\le L\leq \infty$ and $K(\pm 1)=0$.
}

In the future, we denote by $\kappa_j$,  $j=0,1,2,3$, the absolute $j$th moment of the distribution with density $K(t)$, i.e., $\kappa_j= \int\nolimits_{-1}^1|u|^jK(u)du$.
Put
$K_{h}(t)=h^{-1} K(h^{-1}t)$. It is clear that  $K_{h}(s)$ is a probability density with support lying in  $[-h,h]$.
We need also the notation
$$\|K\|^2=\int\limits_{-1}^1K^2(u)du, \qquad\kappa_j(\alpha)=\int\limits_{-1}^{\alpha}t^jK(t)dt,\quad \alpha\in[0,1], \quad j=0,1,2,3.$$

\begin{rem}

We emphasize that assumption $(D)$ includes a fixed design situation.
We consider the segment $[0,1]$ as an area of design change solely for the sake of simplicity of exposition of the approach.
In the general case, instead of the segment $ [0,1] $, one can consider an arbitrary Jordan measurable subset of $\mathbb {R}$.
\end{rem}

Further, we denote by
$z_{n:1}\leq \ldots\leq z_{n:n}$ the order statistics constructed by the sample
$\{z_i;\,i=1,\ldots, n\}$.
Put
 $$z_{n:0}:=0,\quad z_{n:n+1}:=1,\quad  \Delta z_{ni}:=z_{n:i}-z_{n:i-1}, \quad i=1,\ldots,n+1.$$
 The response variable $X_l$ and the random error $\varepsilon_{l}$ from  (\ref{1}) corresponding to the order statistic $z_{n:i}$,
will be denoted by $X_{ni}$ and $\varepsilon_{ni}$, respectively. It is easy to see that the  new errors
 $\{\varepsilon_{ni};\,i=1,\ldots, n\}$ satisfy  condition $(E)$ as well.
Next, by $O_p(\eta_n)$ we denote a random variable
$\zeta_n$ such that, for all $M>0$, one has
   $$\limsup\limits_{n\to\infty}{\mathbb P}(|\zeta_n|/\eta_n>M)\le \beta(M),
   $$
  where  $\lim_{M\to\infty}\beta(M)=0$ and $\{\eta_n\}$ are positive  (maybe random or not) variables and the function $\beta(M)$
  that may depend on the kernel $K$ and  $\sigma^2$.
We agree that, throughout what follows, all limits, unless otherwise stated, are taken for $n\to\infty$.

Let us introduce one more constraint, which is the crucial condition of the paper (in particular, the only condition on design points that guarantees the existence of a uniformly consistent estimator; see also the comments at the end of the section).

 $({\bf D}_0)$
{\it The following limit relation holds$:$
 $\delta_n:=\max\limits_{1\leq i\le n+1}\Delta z_{ni}\stackrel{p}{\to} 0$.
 }

Finally, for any $ h \in (0,1) $,
we introduce into consideration the following class of estimators for the regression function $f$:
\begin{eqnarray}\label{est1}
\widehat f_{n,h}(t):=
I(\delta_n\le c_*h)\sum_{i=1}^n\frac{w_{n2}(t)-(t-z_{n:i})w_{n1}(t)}{w_{n0}(t)w_{n2}(t)-w_{n1}^2(t)}X_{ni}K_{h}(t-z_{n:i})\Delta z_{ni},
\end{eqnarray}
where $I(\cdot)$ is the indicator function,
\begin{equation}\label{2}
c_*\equiv c_*(K):=\frac{\kappa_2-\kappa^2_1}{96
L(6L+\kappa_2+\kappa_1/2)}<\frac{1}{864L};
\end{equation}
hereinafter, we use the notation
$$w_{nj}(t):=\sum_{i=1}^n(t-z_{n:i})^jK_{h}(t-z_{n:i})\Delta z_{ni}, \quad j=0,1,2,3.
$$
\begin{rem}
It is easy to see that the difference $\kappa_2-\kappa^2_1$ is the variance of a non-degenerate distribution, thus this is strictly positive.
\end{rem}

\begin{rem}
It is easy to verify that kernel estimator (\ref{est1}), without the indicator factor, is the first coordinate of the two-dimensional estimate of the weighted least-squares method, i.e., of the two-dimensional point $ (a^*, b^*) $ at which
the following minimum is attained:
\begin{eqnarray}\label{est2}
\min\limits_{a,b}
\sum_{i=1}^n\left(X_{ni}-\big(a+b(t-z_{n:i})\big)\right)^2K_{h}(t-z_{n:i})\Delta z_{ni}.
\end{eqnarray}
Thus, the proposed class of estimators in a certain sense (in fact, by construction) is close to the classical local linear kernel estimators, but in the weighted least squares method (\ref{est2}) we use slightly different weights.
 \end{rem}

\begin{rem}
In the case when there are multiple design points, some spacings $\Delta z_{ni}$ vanish, and we lose some of the sample information in the estimator
 (\ref{est1}).
In this case, it is proposed, before using the estimator (\ref{est1}), to slightly reduce the sample by replacing the observations $X_i$ with the same points $z_i$ with their sample mean and leaving only one design point out of multiples in the new sample.
In this case, the averaged observations will have less noise. So, despite the smaller size of the new sample, we do not lose the information contained in the original sample.
\end{rem}

Let us further agree to denote by $C_j$, $j\geq 1$, absolute positive constants, and by $C_j^*$, positive constants depending only on the kernel $K$.
The main result of this section is as follows.

\begin{theor}\label{theor-1} Let conditions $(D)$, $(E)$, and $(K)$ be satisfied. Then, for any fixed $h\in(0,1/2)$,
with probability $1$ it is satisfied
\begin{eqnarray}\label{5}
\sup_{t\in [0,1]}|\widehat f_{n,t}(t)-f(t)|\le C_1^*\omega_f(h)+\zeta_n(h),
\end{eqnarray}
where $\omega_f(h):=\sup\limits_{u,v\in [0,1]: |u-v|\le
h}|f(u)-f(v)|$ and
the random variable $\zeta_n(h)$ meets the relation
\begin{eqnarray}\label{6}
{\mathbb P}\left( \zeta_n(h)>y,\, \delta_n\le c_*h\right)\le C_2^*\sigma^2\frac{{\mathbb E}\delta_n}{h^2y^2},
\end{eqnarray}
with the constant $c_*$ from $(\ref{2})$.
\end{theor}
\begin{rem}
{\rm As follows from the proof of Theorem 1, the constants $C_1^*$ and $C_2^*$ have the following structure:
$$
C_1^*=C_1\frac{L^2}{\kappa_2-\kappa_1^2},\qquad C_2^*=C_2\frac{L^4}{(\kappa_2-\kappa_1^2)^2}.
$$
}
\end{rem}
\begin{rem}\label{zam-0}
{\rm
Since $\delta_n\le 1$, then under condition $(D_0)$ the limit relation ${\mathbb E}\delta_n\to 0$ holds.
Therefore, taking into account Theorem \ref{theor-1}, we can assert that
$\zeta_n(h)=O_p(h^{-1}({\mathbb E}\delta_n)^{1/2})$.
Thus, the bandwidth $h$
can be determined, for example, by the relation
\begin{eqnarray}\label{7}
h_n=\sup\left\{h>0:\,{\mathbb P}\left(\omega_f(h)\ge h^{-1}({\mathbb E}\delta_n)^{1/2}\right)\le h^{-1}({\mathbb E}\delta_n)^{1/2}\right\}.
\end{eqnarray}
It is easy to see that, when $(D_0)$ is satisfied, the limit relations
$h_n \to 0$ and $h^{-1}_n({\mathbb E} \delta_n)^{1/2} \to 0 $ hold.
In fact, the value of $h_n$
equalizes in $h$ the order of smallness in probability of both terms on the right-hand side of the relation (\ref{5}).
Note also that, for nonrandom $f$, one can choose $h \equiv
{h}_n$ as a solution to the equation
\begin{eqnarray}\label{8}
h^{-1}({\mathbb E}\delta_n)^{1/2}=\omega_f(h).
\end{eqnarray}
It is clear that this solution tends to zero as $n$ grows.

The relations (\ref{7}) and (\ref{8}) allow us to obtain the order of smallness of the optimal bandwidth $h$, but not the optimal value of
$h$. In practice, $h$ can be chosen, for example, by so-called cross-validation.
 \hfill$\square$
}
\end{rem}

From Theorem \ref{theor-1} and Remark \ref{zam-0} it is easy to obtain 
the following corollary.
\begin{cor}
 Let the conditions $(D)$, $(D_0)$, $(K)$, and $(E)$ be satisfied, the regression function $f(t)$
 is nonrandom, and
 $\cal C $ is an arbitrary subset of equicontinuous functions in $C[0,1]$ $($for example, a precompact set$)$. Then
$$
\gamma_n({\cal C})=\sup_{f\in \cal C}\sup_{t\in [0,1]}| \widehat f_{n,\tilde h_n}(t)-f(t)|\stackrel{p}\to 0,
$$
where $\tilde h_n$ is defined by equation $ (\ref{8}) $, in which the modulus of continuity $ \omega_{f}(h) $ is replaced with the universal modulus
$\omega_{\cal C}(h)=\sup\nolimits_{f\in \cal C}\omega_f(h)$.
Moreover, the asymptotic relation $\gamma_n({\cal C})=O_p(\omega_{\cal C}(\tilde h_n))$ hold.
\end{cor}

 \begin{rem}
{\rm It is easy to see that for a nonrandom $ f (t) $ the modulus of continuity in (\ref{8}) can be replaced by one or another upper bound for $\omega _{\cal C}(h)$, obtaining the corresponding upper bound for $\gamma_n({\cal C})$.
Consider the case
${\mathbb E}\delta_n=O(1/n)$.
If $ \cal C $ consists of functions $f(t)$ satisfying the H\"{o}lder condition with exponent $\alpha\in (0,1]$ and a universal constant then
$\tilde h_n = O \left(n^{-\frac{1}{2 (1+ \ alpha)}}\right)$ and
$\omega_{\cal C}(\tilde h_n)=
O\left (n^{-\frac{\alpha}{2(1+\alpha)}}\right )$.
In particular, if the functions from ${\cal C}$ satisfy the Lipschitz condition ($\alpha = 1$) with a universal constant then
$\gamma_n({\cal C})=O_p(n^{-1/4})$.
}
\end{rem}

From Theorem \ref{theor-1} and Remark \ref{zam-0} we obtain 
the following corollary.

\begin{cor}
Let the conditions $ (D) $, $ (D_0) $, $ (K) $, and $ (E) $ be satisfied and let the modulus of continuity $ \omega_{f}(h) $ of the random regression function
$f(t)$ with probability $ 1 $ admit the upper bound $ \omega_{f}(h)\le\zeta d(h)$, where $ \zeta> 0 $ is a random variable and $ d (h) $ is a positive continuous nonrandom function such that
 $d(h)\to 0$ as $h\to 0$. Then
\begin{eqnarray}\label{optim3}
\sup_{t\in [0,1]}| \widehat f_{n,\hat h(n)}(t)-f(t)|\stackrel{p}\to 0,
\end{eqnarray}
where the value $\hat h_n$ is defined in $(\ref{8})$ after replacement
 $d(h)$.
\end{cor}

Let us discuss in more detail condition $(D_0)$.
 Obviously,  condition $(D_0)$ is satisfied for any nonrandom regular design (this is the case of nonidentically distributed $\{z_i\}$
depending on $n$). If $ \{z_i\} $ are independent and identically distributed and the interval $ [0,1] $ is the support of distribution of $ z_1 $
then  condition $ (D_0) $ is also satisfied. In particular, if the distribution density of $ z_1 $ is separated from zero on $ [0,1] $ then
 $ \delta_n = O\left(\log n/n\right) $ holds (see details in \cite{2020}).
 If $ \{z_i; \, i \ge 1 \} $ is a stationary sequence with a marginal distribution with the support $[0,1]$, satisfying an $ \alpha $-mixing condition then
 condition $ (D_0) $  is also satisfied (see Remark  \ref{rem-4} below).
Note that the dependence of the random variables $\{z_i\}$ satisfying  condition $ (D_0) $ can be much stronger,
which is illustrated in the following example.
 \begin{ex}\label{ex-1}
 {\rm
Let the sequence of random variables $ \{z_i; \, i\ge 1 \} $ be defined by the relation
\begin{equation}\label{11}
{ z}_{i}=\nu_{i}{u}_{i}^{l}+(1-\nu_{i}){u}_{i}^{r},
\end{equation}
where $\{{u}_{i}^l\}$ and  $\{{u}_{i}^r\}$ are independent and uniformly distributed on $ [0,1/2] $ and $ [1/2,1] $, respectively, the sequence
 $ \{\nu_i \} $ does not depend on $ \{{u}_{i}^l \} $, $ \{{u}_{i}^r \} $ and consists of Bernoulli random variables with success  probability $1/2 $, i.e., the distribution of random variables $ {z}_i $ is an equilibrium mixture of two uniform distributions on the corresponding intervals.
The dependence between the random variables $ \nu_i $ for any natural number $ i $ is defined by the equalities
$ \nu_{2i-1} = \nu_1 $ and $ \nu_{2i} = 1- \nu_1 $.
In this case, the random variables $\{z_i;\,i\ge 1\}$ in (\ref{11}) form a stationary sequence of random variables uniformly distributed on the segment $ [0,1] $, satisfying condition $ (D_0) $.
On the other hand,
for all natural numbers $m$ and $n$,
\begin{eqnarray*}
{\mathbb P}(z_{2m}\le 1/2,\,z_{2n-1}\le 1/2)=0.
\end{eqnarray*}
Thus, all the known conditions for the weak dependence of random variables (in particular, the mixing conditions) are not satisfied here.

According to the scheme of this example, it is possible to construct various sequences of dependent random variables uniformly distributed on $ [0,1] $ by
choosing sequences of Bernoulli switches with the conditions $ \nu_{j_k} = 1 $ and $ \nu_{l_k} = 0 $ for infinite numbers of
indices $ \{j_k\} $ and $ \{l_k\} $. In which case, condition $ (D_0) $ will also be satisfied, but the corresponding sequence
$ \{z_i\} $
(not necessarily stationary) may not even satisfy the strong law of large numbers.
For example, this is the case
when $ \nu_j = 1-\nu_1 $ for $ j = 2^{2k-1}, \ldots, 2^{2k}-1 $, and $ \nu_j = \nu_1 $ for $ j = 2^{2k} , \ldots, 2^{2k + 1}-1 $, where $ k = 1,2, \ldots $ (i.e., we randomly choose one of the two segments $ [0,1/2] $ and $ [1/2,1] $,
into which we randomly throw the first point, and then alternate the selection of one of the two segments
by the following numbers of elements of the sequence: $1$, $ 2$, $2^2$, $2^3$,  etc.).
Indeed, we can introduce the notation $n_k = 2^{2k}-1$, $\tilde n_k = 2^{2k + 1}-1 $, $S_m = \sum\nolimits_{i = 1}^m z_{i}$ and note that, for all elementary events from the event $\{\nu_1 = 1\}$, one has
\begin{equation*}
\frac{S_{n_k}}{n_k}=\frac{1}{n_k}\sum\limits_{i\in N_{1,k}}u_{i}^l+\frac{1}{n_k}\sum\limits_{i\in N_{2,k}}u_{i}^r,
\end{equation*}
where   $N_{1,k}$ and $N_{2,k}$  are the sets of indices, for which the observations $\{{z}_i, i\leq n_k\}$ lie in the  intervals  $[0,1/2]$ or $[1/2,1]$, respectively. It is easy to see that $\#(N_{1,k})=n_k/3$ and   $\#(N_{2,k})=2\#(N_{1,k})$. Hence,    ${S_{n_k}}/{n_k}\to{7}/{12}$ almost surely as  $k\to\infty$ due to the strong law of large numbers for the sequences  $\{u_{i}^l\}$ and  $\{u_{i}^r\}$. On the other hand, as $k\to\infty$, for all elementary events  from $\{\nu_1=1\}$  one has
 \begin{equation}\label{10}
\frac{S_{\tilde n_k}}{\tilde n_k}=\frac{1}{\tilde n_k}\sum\limits_{i\in \tilde N_{1,k}}u_{i}^l+\frac{1}{\tilde n_k}\sum\limits_{i\in \tilde N_{2,k}}u_{i}^r\to \frac{5}{12},
\end{equation}
where  $\tilde N_{1,k}$ and $\tilde N_{2,k}$ are the sets of indices, for which the observations $\{{z}_i, i\leq \tilde n_k\}$ lie in the intervals  $[0,1/2]$ or  $[1/2,1]$, respectively. Proving the convergence in  (\ref{10}), we took into account that $\#(\tilde N_{1,k})=(2^{2k+2}-1)/3$ and $\#(\tilde N_{2,k})=2n_k/3$, i.e.,  $\#(\tilde N_{1,k})=2\#(\tilde N_{2,k})+1$.

Similar arguments are valid for all elementary events from $\{\nu_1=0\}$.
$\hfill\Box$
}
\end{ex}

\begin{rem} \label{rem-4}
In the case of i.i.d. random variables  $\{z_i\}$,   condition $(D_0)$ will be fulfilled if, for all $\delta\in (0,1)$,
\begin{equation}\label{izmel2}
p_n(\delta)\equiv\sup_{|\Delta|=\delta}{\mathbb P}\Big(\bigcap\limits_{i\le n}\{z_i\notin \Delta\}\Big)\to 0,
\end{equation}
where the supremum is taken over all intervals $\Delta\subset [0,1]$ of length $\delta$.
Indeed, for any natural $N>1$, we divide the interval $[0,1]$ into $N$ subintervals $\Delta_k$, $k=1,\ldots,N$, of length  $1/N$.
Then one has
$${\mathbb P}\Big(\max\limits_{1\leq i\le n+1}\Delta z_{ni}>\frac{2}{N}\Big)\leq \sum\limits_{k=1}^N{\mathbb P}\Big(\bigcap\limits_{i\le n}\{z_i\notin \Delta_k\}\Big)\leq N \max\limits_k{\mathbb P}\Big(\bigcap\limits_{i\le n}\{z_i\notin \Delta_k\}\Big)\leq Np_n(1/N),$$
since the event $\big\{\max\nolimits_{1\leq i\le n+1}\Delta z_{ni}>2/N\big\}$ implies the existence of an interval $\Delta_k$ of length~$1/N$ that does not contain any points from the collection $\{z_i\}$.
Thereby, condition (\ref{izmel2}) implies the limit relation $ \max\nolimits_{i \le n + 1} \Delta z_{ni} \stackrel{p}\to 0 $, which is equivalent to convergence with probability $1$ due to the monotonicity of the sequence $\max\nolimits_{i \le n + 1} \Delta z_{ni}$.
In particular, if $ \{z_i\} $ are independent then $ p_n(\delta) = e^{-c(\delta)n}$ and $c(\delta)> 0$, i.e., as $ n\to\infty $, the finite collection $\{z_i\}$ with probability $ 1 $ form a refining partition of the finite segment $ [0,1] $. It is easy to show that if $\{z_i;\,i \ge 1 \} $ is a stationary sequence satisfying an $\alpha$-mixing condition and having a marginal distribution with the support $[0,1]$ then (\ref{izmel2}) will be valid.
  \hfill$\square$
\end{rem}

\section{Estimating the mean function
of a stochastic
process}

Consider the following statement of the problem of estimating the expectation of an almost surely continuous stochastic  process $f(t)$.
There are $ N $ independent copies of the regression equation (\ref {1}):

 \begin{equation}  \label{50}
 X_{i,j}=f_j(z_{i,j})+\varepsilon_{i,j},
 \qquad i=1,\ldots,n,\,\,\,\,j=1,\ldots,N,
 \end{equation}
 where $f(t), f_1(t),\ldots, f_N(t)$, $t\in [0,1]$, are
 independent identically distributed  almost surely continuous unknown random processes, the set $ \{\varepsilon_{i, j}; \,i = 1, \ldots, n \} $ satisfies  condition $ (E) $ for any $ j $, the set $ \{z_{i, j}; \, i = 1, \ldots, n \} $ meets conditions $ (D) $ and $ (D_0) $ for any $ j $ (here and below the index $ j $ for the considered random variables means the number of copy of Model (\ref{1})).
In particular, under the assumption that condition $ (K) $ is valid, by $ \widehat f_{n, h, j} (t) $, $ j = 1, \ldots, N, $ we denote the estimator given by the relation (\ref{est1}) when replacing the values from (\ref{1}) with the corresponding characteristics from (\ref{50}). Finally,
an estimator for the mean-function is determined by the equality
\begin{equation}\label{LLN}
\overline{\widehat f_{N,n,h}}(t)=\frac{1}{N}\sum\limits_{j=1}^N\widehat f_{n,h,j}(t).
\end{equation}

As a consequence of Theorem \ref{theor-1}, we obtain the following assertion.
\begin{theor}\label{theor-2}
Let Model  $(\ref{50})$ satisfy the above-mentioned conditions and, moreover,
\begin{equation}\label{51}
{\mathbb E}\sup_{t\in [0,1]}|f(t)|<\infty,
\end{equation}
while the sequences  $h\equiv h_n\to 0$ and  $N\equiv N_n\to \infty$ meet the restrictions
\begin{equation}\label{52}
h^{-2}{\mathbb E}\delta_n\to 0\,\,\,\mbox{and}\,\,\,N{\mathbb P}(\delta_n>c_*h)\to 0.
\end{equation}
Then
\begin{equation}\label{53}
\sup\limits_{t\in [0,1]}\left|\overline{\widehat f_{N,n,h}}(t)-{\mathbb E}f(t)\right|\stackrel{p}\to 0.
\end{equation}
\end{theor}

\begin{rem}\label{za3}
 {\rm If condition (\ref {51}) is replaced with a slightly  stronger constraint
 $${\mathbb E}\sup\nolimits_{t\in [0,1]}f^{2}(t)<\infty$$
then, under conditions similar to (\ref{52}), one can prove the uniform consistency of the estimator
$$\widehat M_{N,n,h}(t_1,t_2)=\frac{1}{N}\sum\limits_{j=1}^N\widehat f_{n,h,j}(t_1)\widehat f_{n,h,j}(t_2),\,\,\,\,t_1,t_2\in [0,1],
$$
for the unknown mixed second moment ${\mathbb E}f(t_1)f(t_2)$
where $h\equiv h_n$ and  $N\equiv N_n$ satisfy (\ref{52}). The arguments in proving this fact are quite similar to those in proving Theorem \ref{theor-2} and they are omitted.
In other words, under the above-mentioned restrictions, the estimator
$${\widehat{ \rm Cov}}_{N,n,h}(t_1,t_2)= \widehat M_{N,n,h}(t_1,t_2)-\overline{\widehat f_{N,n,h}}(t_1)\overline{\widehat f_{N,n,h}}(t_2)$$ uniformly consistent for the covariance of the random regression function $f(t)$.
}
\end{rem}

\begin{rem}
The problem of estimating the mean and covariance functions plays a fundamental role in the so-called functional data analysis (see, for example, \cite{HE-2015}, \cite{KR-2017}, \cite{2010-Li}, \cite{Mu-2005}). The property of uniform consistency of certain estimates of the mean function, which is important in the context of the problem under consideration, was considered, for example, in \cite{HE-2015}, \cite{2010-Li}, \cite{YMW-2005}, \cite{ZW-2016}, \cite{ZLL-2018}.
  For a random design, as a rule, it is assumed that all its elements are independent identically distributed random variables (see, for example, \cite{CY-2011}, \cite{2006-Ha}, \cite{2010-Li},   \cite{WZ-2006}--\cite{ZLL-2018}).
  In the case where the design is deterministic, certain regularity conditions discussed above in Introduction are usually used.
 Moreover, in the problem of estimating the mean function, it is customary to subdivide design elements into certain types depending on the density of filling with the
design points the regression function domain.
The literature focuses on two types of data:
 or the design is in some sense ``sparse'' (for example, the number of design elements in each series is uniformly limited  \cite{CY-2011}, \cite{2006-Ha}, \cite{2010-Li}, \cite{WZ-2006},  \cite{ZLL-2018}), or the design is somewhat ``dense'' (the number of elements in each series grows with the number of series \cite{CWLY-2016}, \cite{2010-Li}, \cite{WZ-2006}, \cite{ZC-2007},  \cite{ZLL-2018}). Theorem \ref{theor-2} considers the second of the specified types of design under condition $(D_0)$ in each of the independent series.
Note that our formulation of the problem of estimating the mean function also includes the situation of a general deterministic design.

Note that the methodologies for estimating the mean function used for dense or sparse data are often different
(see, for example, \cite{Mu-2005}, \cite{WCM-2016}). In the situation of a growing number of observations in each series, it is natural to preliminarily estimate
trajectories of a random regression function in each series, and then average over all series (e.g., see \cite{CY-2011}, \cite{2006-Ha},
  \cite{ZC-2007}). This is exactly what we do in (\ref{LLN}) following this conventional approach. \hfill$\square$

\end{rem}

\section{Comparison with some known approaches
}

In \cite{2020},  under the conditions of the present paper, the following estimators were studied:
\begin{equation}\label{est4}
f^*_{n,h}(t)=\frac{\sum_{i=1}^nX_{ni}K_{h}(t-z_{n:i})\Delta z_{ni}}{\sum_{i=1}^nK_{h}(t-z_{n:i})\Delta z_{ni}}\equiv \frac{\sum_{i=1}^nX_{ni}K_{h}(t-z_{n:i})\Delta z_{ni}}{w_{n0}(t)}.
\end{equation}
Notice that
\begin{equation}\label{est5}
f^*_{n,h}(t)\equiv {\rm arg}\min\limits_{a}
\sum\limits^n_{i=1}(X_{ni}-a)^2K_{h}(t-z_{n:i})\Delta z_{ni}.
\end{equation}

It is interesting to compare the new estimators $\widehat f_{n,h}(t)$ with the estimators $f_{n,h}^*(t)$ from \cite{2020} as well as with other estimators (for example, the Nadaraya--Watson estimators $\widehat f_{NW}(t)$ and classical local linear estimators $\widehat f_{LL}(t)$).
Throughout this section, we assume that conditions $ (D) $, $ (K) $, and $ (E) $ are satisfied and the regression function $ f(t) $
 is nonrandom. Moreover,
  we need the following constraint.

 ${({\bf IID})}$ {\it The regression function $ f (t) $ in Model {\rm(\ref{1})}
 twice continuously differentiable, the errors $ \{\varepsilon_i\} $ are independent, identically distributed, centered, and independent of the design
 $\{z_i\}$,
whose elements are independent and identically distributed. In addition, the distribution function of the random variable $ z_1 $ has a strictly positive density $ p(t) $ continuously differentiable on $ (0,1) $.
}

Such severe restrictions on the parameters of the regression model are explained both by
   problems in calculating the asymptotic representation for the variances of the estimators $ \widehat f_{n,h}(t) $ and $ f_{n,h}^*(t) $ as well as
   by properties of the Nadaraya--Watson estimators, which are very sensitive to the nature of the correlation of design elements.

For any statistical estimator $ \tilde f_n(t) $ of the regression function $ f(t) $, we will use the notation $ {\rm Bias} \tilde f_n(t) $ for its bias, i.e.,
$
{\rm Bias}\tilde f_n(t):={\mathbb E}\tilde f_n(t)-f(t).
$
Put $\overline f =\sup\nolimits_{t\in [0,1]}|f(t)|$
and for $j=0,1,2,3$, introduce the notation
\begin{eqnarray}\label{int}
w_{j}(t)=\int\limits_0^1(t-z)^jK_{h}(t-z)dz=\int\limits_{z\in [0,1]: |t-z|\leq h}(t-z)^jK_{h}(t-z)dz,\quad t\in [0,1].
\end{eqnarray}

The following asymptotic representation
for the bias and variance of the estimator $ f_{n,h}^*(t) $ was obtained in \cite{2020}.

        \begin{prop} \label{le2vve}
Let condition ${ (IID)}$ be fulfilled and $\inf_{t\in[0,1]}p(t)>0$. If  $n\to\infty$ and  $h\to 0$ so that  $(\log n)^{-1}h\sqrt n\to\infty$, $h^{-2}{\mathbb E}\delta_n\to 0$, and $h^{-3}{\mathbb E}\delta_n^2\to 0$ then,
for any  $t\in (0,1)$, the following asymptotic relations are valid$:$
$$
{\rm Bias}f_{n,h}^*(t)=\frac{h^2\kappa_2}{2}f''(t)+o(h^2),\qquad
{\mathbb Var}f_{n,h}^*(t)\sim \frac{2\sigma^2}{h np(t)}\|K\|^2.
$$
 \end{prop}

Note that the first statement concerning the asymptotic behavior of the bias in Proposi\-tion~\ref{le2vve} was actually proved for arbitrarily dependent
design elements when condition $ (D_0) $ is met. The following two propositions and corollaries are also obtained without any assumptions about
correlation of design elements, only conditional centering and conditional orthogonality of the errors from condition $ (E) $ are used.

 \begin{prop} \label{predl-4}
 Let $ h<1/2$.
Then, for any fixed $t\in [h,1-h]$,
$$ {\rm Bias}\widehat f_{n,h}(t)=  {\rm Bias} f_{n,h}^*(t)+\gamma_{n,h}(t),\quad
{\mathbb Var}\widehat f_{n,h}(t)={\mathbb Var}f_{n,h}^*(t)+\rho_{n,h}(t),
$$
where $$|\gamma_{n,h}(t)|\le C_3^*\overline f h^{-1}{{\mathbb E}\delta_n},\quad
|\rho_{n,h}(t)|\le C_4^*\big(\sigma^2+\overline f^2 \big)h^{-1}{{\mathbb E}\delta_n}.$$
 \end{prop}

 \begin{prop} \label{predl-5}
Let the regression function $ f(t) $ be twice continuously differentiable. Then, for any fixed $ t\in (0,1) $,
\begin{equation}\label{hatf}
{\rm Bias}\widehat f_{n,h}(t)=\frac{f''(t)}{2}B_{0}(t)+O({\mathbb E}\delta_n/h)+o(h^2),
\end{equation}
where
\begin{equation}\label{B_0}
B_0(t)=\frac{w^2_{2}(t)-w_{3}(t)w_{1}(t)}{w_{0}(t)w_{2}(t)-w^2_{1}(t)}.
\end{equation}

Moreover,
\begin{equation}\label{starf}
{\rm Bias} f_{n,h}^*(t)=-f'(t)\frac{w_{1}(t)}{w_{0}(t)}+\frac{f''(t)}{2}\frac{w_{2}(t)}{w_{0}(t)}
+ O({\mathbb E}\delta_n)+o(h^2),
\end{equation}
besides, the error terms $o(h^2)$ and $O(\cdot)$ in $(\ref{hatf})$ and $(\ref{starf})$ are uniform in $t$.
\end{prop}

\begin{cor}\label{equiv}
 {\it Let the regression function $ f(t) $ be twice continuously differentiable,
$h\to 0$, and  $h^{-3}{\mathbb E}\delta_n\to 0$. Then, for each fixed $ t\in (0,1) $ such that
 $f''(t)\neq 0$, the following asymptotic relations are valid$:$}
$${\rm Bias}\widehat f_{n,h}(t)\sim {\rm Bias} f_{n,h}^*(t)\sim \frac{f''(t)}{2}\kappa_2h^2.
$$
\end{cor}

 \begin{cor}\label{atzero}
Suppose that, under the conditions of the previous corollary, $ f $ has nonzero first and second derivatives in a neighborhood of zero. Then for any fixed
positive $ \alpha <1 $ such that $ \kappa_1 (\alpha) <0 $, the following asymptotic relations hold$:$
$${\rm Bias}\widehat f_{n,h}(\alpha h)\sim \frac{1}{2}h^2D(\alpha)f''(0+),\qquad
{\rm Bias} f_{n,h}^*(\alpha h)\sim -h\frac{\kappa_1(\alpha)}{\kappa_0(\alpha)}f'(0+),$$
where
$$D(\alpha)=\frac{\kappa^2_2(\alpha)-\kappa_3(\alpha)\kappa_1(\alpha)}{\kappa_0(\alpha)\kappa_2(\alpha)-\kappa^2_1(\alpha)}.$$
\end{cor}
Note that, due to the Cauchy--Bunyakovsky inequality and the properties of the density $ K (\cdot) $, the strict inequality $ \kappa_0 (\alpha) \kappa_2 (\alpha) - \kappa^2_1(\alpha)> 0 $ holds for any $ \alpha \in [0,1] $.

\begin{rem} Similar relations take place in a neighborhood of the right boundary of the segment $ [0,1] $, when $ t = 1- \alpha h $ for any $ \alpha \le 1 $. In this case, in the above asymptotics, one simply needs to replace the right-hand derivatives at zero by analogous (non-zero) left-hand derivatives at point 1, and instead of the quantities
$ \kappa_j(\alpha) $ must be substituted $ \tilde \kappa_j (\alpha) = \int \nolimits_{-\alpha}^1v^iK(v)dv = (-1)^j \kappa_j(\alpha) $. In this case, the coefficient $ D (\alpha) $ will not change, and the corresponding coefficient on the right-hand side of the second asymptotics will only change its sign.
\end{rem}

Thus, the qualitative difference between the estimators $ f_{n,h}^*(t)$ and $ \widehat f_{n,h}(t) $ is observed only in neighborhoods of the boundary points $ 0 $ and $ 1 $: For the estimator $ f_{ n,h}^*(t)$, in the $ h $-neighborhoods of the indicated points, the order of smallness of the bias is $h$, and for $ \widehat f_{n,h}(t) $ this order is $h^2$. Such a connection between the estimators (\ref{est1}) and (\ref{est4}) seems to be quite natural in view of the relations (\ref{est2}) and (\ref{est5}),
and the known relationship at the boundary points between  Nadaraya--Watson estimators $ \hat f_{NW}(t) $ and  locally linear estimators $\widehat f_{LL}(t)$.

\begin{rem} If condition $ {(IID)} $ is satisfied then, for the bias and variance of estimators $ \widehat f_{NW}(t) $ and $ \widehat f_{LL}(t) $,
  the following asymptotic representations are well known (see, for example, \cite{1996-Fa}), which are valid for any $ t\in (0,1) $  under broad conditions on the parameters of the model under consideration:
\begin{eqnarray*}
{\rm Bias}\widehat f_{NW}(t)=\frac{h^2\kappa_2}{2p(t)}\left( f''(t)p(t) + 2f'(t)p'(t)\right)+o(h^2), \quad
{\mathbb Var}\widehat f_{NW}(t)\sim \frac{\sigma^2}{h np(t)}\|K\|^2,\\
{\rm Bias} \widehat f_{LL}(t)=\frac{h^2\kappa_2}{2}f''(t)+o(h^2),\qquad
{\mathbb Var}\widehat  f_{LL}(t)\sim \frac{\sigma^2}{h np(t)}\|K\|^2.
\end{eqnarray*}

The above asymptotic representations show that if the assumptionss $ {(IID)}$ are valid then
the variance of the Nadaraya--Watson estimator $ \hat f_{NW}(t) $ and the locally linear  estimator $ \hat f_{LL}(t) $ under broad conditions is asymptotically half the variance of the estimators $ f^*_{n,h}(t) $ and $ \hat f_{n,h}(t) $, respectively.
But the mean-square error of any estimator is equal to the sum
of the  variance and squared bias, which for the compared estimators
is asymptotically determined by the quantities $ f''(t)p(t) + 2f'(t)p'(t) $ or $ f''(t)p(t) $, respectively. In other words, if the standard deviation $\sigma$ of the errors is not very large and
\begin{equation}\label{ineqcomp}
\left|f''(t){p(t)}+2f'(t){p'(t)}\right|>
\left|f''(t)p(t)\right|,
\end{equation}
then the estimator $ f^*_{n,h}(t) $ or $ \hat f_{n,h}(t) $ may be more accurate than $ \hat f_{NW}(t) $.
The indicated effect for the estimator $ f^*_{n,h}(t) $ is confirmed by the results of computer simulations in \cite{2020}.

Note also that in order to choose in a certain sense the optimal bandwidth $ h $, the orders of the smallness of the bias and the standard deviation of the estimator are usually equated. In other words, if the assumptions $ {(IID)} $ are fulfilled, for all four types of estimators considered here, we need to solve the equation
$ h^2 \approx (nh)^{-1/2} $. Thus the optimal bandwidth has the standard order
 $h\approx n^{-1/5}$. \hfill$\square$
\end{rem}
 \begin{rem}\label{13}
 Estimators of the form $ \widehat f_{n,h}(t) $ and $ f^*_{n,h}(t) $ given in (\ref{est1}) and (\ref{est4}) can define a little differently, depending on the choice of one or another partition with highlighted points $ \{z_{i};\, i=1,\ldots,n\} $ of the domain of the regression function underlying these estimators.
 For example, using the Voronoi partition of the segment $ [0,1] $,
 an estimator of the form (\ref{est4}) can be given by the equality
\begin{equation}\label{est4+}
\widetilde f^*_{n,h}(t)=\frac{\sum_{i=1}^nX_{ni}K_{h}(t-z_{n:i})
\widetilde\Delta z_{ni}}{\sum_{i=1}^nK_{h}(t-z_{n:i})
\widetilde\Delta z_{ni}},
\end{equation}
where  $\;\widetilde\Delta z_{n1}={\Delta }z_{n1}+{{\Delta }z_{n2}}/{2}$,
$\;\widetilde\Delta z_{nn}={{\Delta }z_{nn}}/{2}+{\Delta }z_{nn+1}$,
$\widetilde\Delta z_{ni}={({\Delta }z_{ni}+{\Delta }z_{ni+1})}/{2}$ for  $i=2,\ldots,n-1$.
Looking through the proofs from \cite{2020} it is easy to see that in this case all properties of the estimator $ \widetilde f_{n,h}^*(t) $ are preserved, except for the asymptotic representation of the variance. Repeating with obvious changes the arguments in proving Proposition \ref{le2vve} in \cite{2020}, we have
 $${\mathbb Var}\widetilde f_{n,h}^*(t)\sim \frac{1.5\sigma^2}{h np(t)}\|K\|^2.$$
 Thus, in the case of independent and identically distributed design points, the asymptotic variance of the estimator can be somewhat reduced by choosing one or another partition.

Similarly, in the definition (\ref{est1}), the estimators $ \widehat f_{n,h}(t) $, the quantities $ \{\Delta z_{ni} \} $ can be replaced by the Voronoi tiling $ \{\widetilde \Delta z_{ni} \} $. It is also worth noting that the indicator factor involved in the determination (\ref{est1}) of the estimator
$ \widehat f_{n,h}(t) $,
does not affect the asymptotic properties of the estimator given in Theorem \ref{theor-1}, and we only needed it to calculate the exact asymptotic behavior of the estimator bias.
\hfill$\square$
 \end{rem}

\section{Simulations
} \label{sec-sim}

In the following computer simulations,
instead of estimator \eqref{est1}, we used the equivalent estimator $\widehat f_{n,h}(t)$ of the weighted least-squares method defined by the relation
\begin{eqnarray}\label{ull1}
(\widehat f_{n,h}(t),\hat b(t))&=&
 {\rm arg}\min\limits_{a,b}
\sum_{i=1}^n\left(X_{ni}-a-b(t-z_{n:i})\right)^2K_{h}(t-z_{n:i})\widetilde\Delta z_{ni},
\end{eqnarray}
where the quantities $\widetilde\Delta z_{ni}$ are defined in (\ref{13}) above. Estimator (\ref{ull1}) differs from  estimator (\ref{est1}) by excluding the indicator
factor and replacing $\Delta z_{ni}$ with $\widetilde\Delta z_{ni}$, which is not essential (see Remark~\ref{13}). Besides, if we had several observations at one design point,
then the observations were replaced by one observation presenting their arithmetic mean (see Remark 4 above).
Although the notation $\widehat f_{n,h}(t)$ in (\ref{ull1}) is somewhat different from the same notation
in \eqref{est1}, we retained the notation $\widehat f_{n,h}(t)$,
which will not lead to ambiguity.

In the simulations below, we will also consider the local constant estimator $\widetilde f^*_{n,h}(t)$ from (\ref{est4+}),
which can be defined by the equality
\begin{equation}\label{ull0}   
\widetilde f^*_{n,h}(t)\equiv {\rm arg}\min\limits_{a}
\sum\limits^n_{i=1}(X_{ni}-a)^2K_{h}(t-z_{n:i})\widetilde\Delta z_{ni}.
\end{equation}
Here we also replace the observations corresponding to one design point 
by their arithmetic mean.

Recall that the Nadaraya-Watson estimator differs from (\ref{ull0}) by
the absence of the factors $\widetilde\Delta z_{ni}$ in the weighting coefficients:
\begin{equation}\label{est4nv}
 \widehat f_{NW}(t)=\frac{\sum_{i=1}^nX_{ni}K_{h}(t-z_{n:i})}{\sum_{i=1}^nK_{h}(t-z_{n:i})}.
\end{equation}
The Nadaraya-Watson estimators are also weighted least-squares estimators:
\begin{equation}\label{est5nv}
\widehat f_{NW}(t)\equiv {\rm arg}\min\limits_{a}
\sum\limits^n_{i=1}(X_{ni}-a)^2K_{h}(t-z_{n:i}).
\end{equation}

In the following examples,  estimators (\ref{ull1}) and (\ref{ull0}), which will be called {\it universal local linear} (ULL) and {\it universal
local constant} (ULC), respectively, will be compared with the estima\-tor of linear regression (LR),
the Nadaraya-Watson (NW) estimator, LOESS of order 1,
 as well as with estimators of generalized additive models (GAM) and of random forest (RF).
 For LOESS estimators, the R {\it loess}() function was used.

It is worth noting that, in the examples below, the best results were obtained by the new estimators
(\ref{ull1}) and (\ref{ull0}), LOESS estimator of order 1, and
the Nadaraya-Watson estimator.

With regard to the simulation examples, the main difference between
estimators (\ref{ull1}) and (\ref{ull0}), and the Nadaraya--Watson and LOESS ones is that
estimators (\ref{ull1}) and (\ref{ull0}) are ``more local''.
This means that if a function $f(z)$ is evaluated on a design interval
$A$ with a ``small'' number of observations adjacent to a design interval $B$
with a ``large'' number of observations, the Nadaraya-Watson and LOESS estimators
will primarily seek to adjust to the ``large'' cluster of observations on
the interval $B$. At the same time,  estimators (\ref{ull1}) and (\ref{ull0}) will equally consider
observations on intervals of equal lengths, regardless of the distribution of design points on the intervals.

In the examples below, for all of the kernel estimators which are the Nadaraya-Watson ones, LOESS,
(\ref{ull1}), and (\ref{ull0}), we used the tricubic kernel
$$K(t)=\frac{70}{81}\max\{0,(1 - |t|^3)^3\}.$$
We chose the tricubic kernel because that kernel is employed in the $R$ function
{\sf loess()} which was
used in the simulations.

The accuracy of the models was estimated with respect to the maximum error and the mean squared error.
In all the examples below, except Example 3, the maximum error was estimated on the uniform grid of 1001 points on the segment $[0,10]$ by the formula
$$\max_{j=1,\dots,1001} |\check f (t_j) - f(t_j)|, $$
where $t_j$ are the grid points of segment $[0,10]$, $t_1=0$, $t_{1001}=10$,
$\check f (t_j)$ are the values of the constructed estimator at the points of the partition grid,
$f (t_j)$ are the true values of the estimated function.
In  Example 3, a grid of 1001 points was taken on the interval from the minimum to the maximum point of the design. That was done in order to
to avoid assessing the quality of extrapolation, since,
in that example, the minimum design point could fall far from 0.

The mean squared error
was calculated for one random splitting of the whole sample into training and validation samples in proportion of $80\%$ to $20\%$, according to the formula
$$\frac{1}{m} \sum_{j=1}^m \left(\check f(z_j) - X_j)\right)^2,$$
where $m$ is the validation sample size, $z_j$ are the validation sample design points,
$X_j$ are the noisy observations of the predicted function in the validation sample,
$\check f$ is the estimate calculated by the training sample.
The splittings into training and validation samples were identical for all models.

For each of the kernel estimators, the parameter $h$ of the kernel $K_h$ was determined using cross-validation minimizing the mean squared error, where the set of observations was partitioned into 10 folds randomly.
The same partitions were taken for all the kernel estimators.

When calculating the root mean square error, the cross-validation for choosing $h$ was carried out on the training set.
To calculate the maximum error, the cross-validation was performed on the whole sample.
For the Nadaraya-Watson models as well as for estimatiors (\ref{ull1}) and (\ref{ull0}),
the parameter $h$ was selected from 20 values located on the logarithmic grid
from $\max\{0.0001, 1.1  \max_i\Delta z_{ni}\}$ to 0.9.
For LOESS, the parameter {\it span} was chosen in the same way from 20 values located on the logarithmic grid from 0.0001 to 0.9.

The simulations also included testing basic statistical learning algorithms: linear regression without regularization,
generalized additive model, and random forest \cite{elstat-2009}. The training of the generalized additive model was carried out using the R library {\it mgcv}.
Thin-plate splines were used, the optimal form of which was selected using generalized cross-validation. Random forest training was done
using the R library {\it randomForest}. The number of trees is chosen to be 1000 based on the out-of-bag error plot for a random forest with five
observations per leaf. The optimal number of observations in a random forest leafs was chosen using 10-fold cross-validation on a logarithmic grid
out of 20 values from 5 to 2000.

In each example, 1000 realizations of different train and vadidation sets were performed, for each of which the errors were calculated.
In each of  train and vadidation sets realizations, 5000 observations were generated.
The results of the calculations are presented below in the boxplots,
where every box represents the median and the 1st and 3rd quartiles.
The plots do not show the results of linear regression, since in the examples, the results appeared to be significantly worse than those of the other models.
The mean squared and maximum errors of estimator \eqref{ull1} were compared with
the errors of LOESS estimator by the paired Wilcoxon test.
The summaries of the errors on the 1000 realizations of different train and vadidation sets
are reported as median (1st quartile, 3rd quartile).

The examples of this section were constructed so that the distribution of design points is ``highly nonuniform''. Potentially, this could demonstrate the advantage of the new estimator \eqref{ull1} over known estimation approaches.

\bigskip

\begin{ex}
{\rm
Let us set the target function
\begin{equation}\label{exam1}
f(z) = (z-5)^2 + 10, \quad 0\le z\le 10
\end{equation}
and let the noise be centered Gaussian with standard deviation $\sigma=2$
(Fig.~\ref{fig1}).
In each realization, we draw 4500
independent design points uniformly distributed on the segment
$z\in [0,5]$,
and
500
independent design points uniformly distributed on the segment $z\in [5,10]$.

\begin{figure}[!ht]
  \centering
  \includegraphics[width=4in]{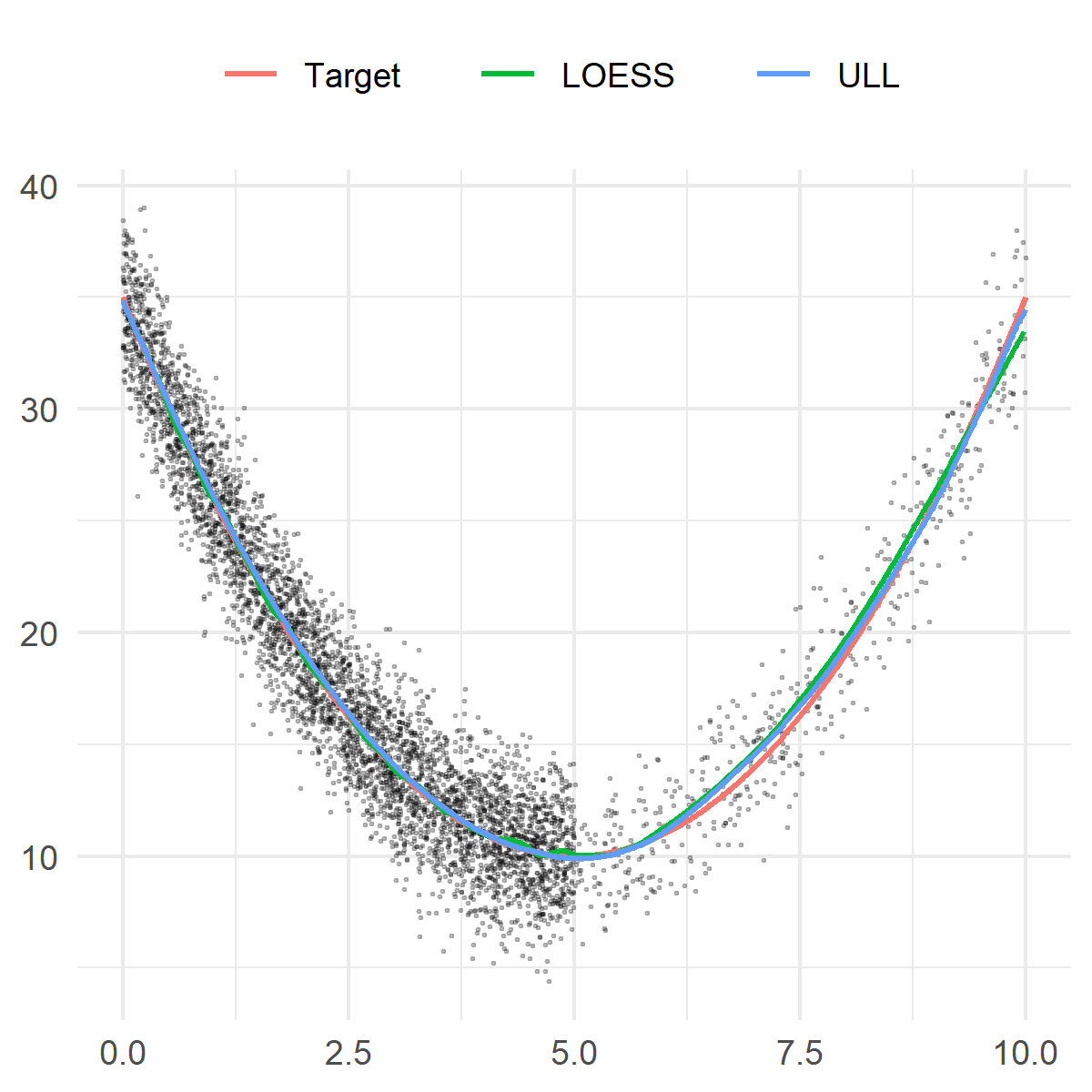}
  \caption{\footnotesize Example 2. Sample observations, target function, and two estimators.}\label{fig1}
\end{figure}

\begin{figure}[!ht]
\centering
  \subfigure{
    \includegraphics[width=.45\textwidth]{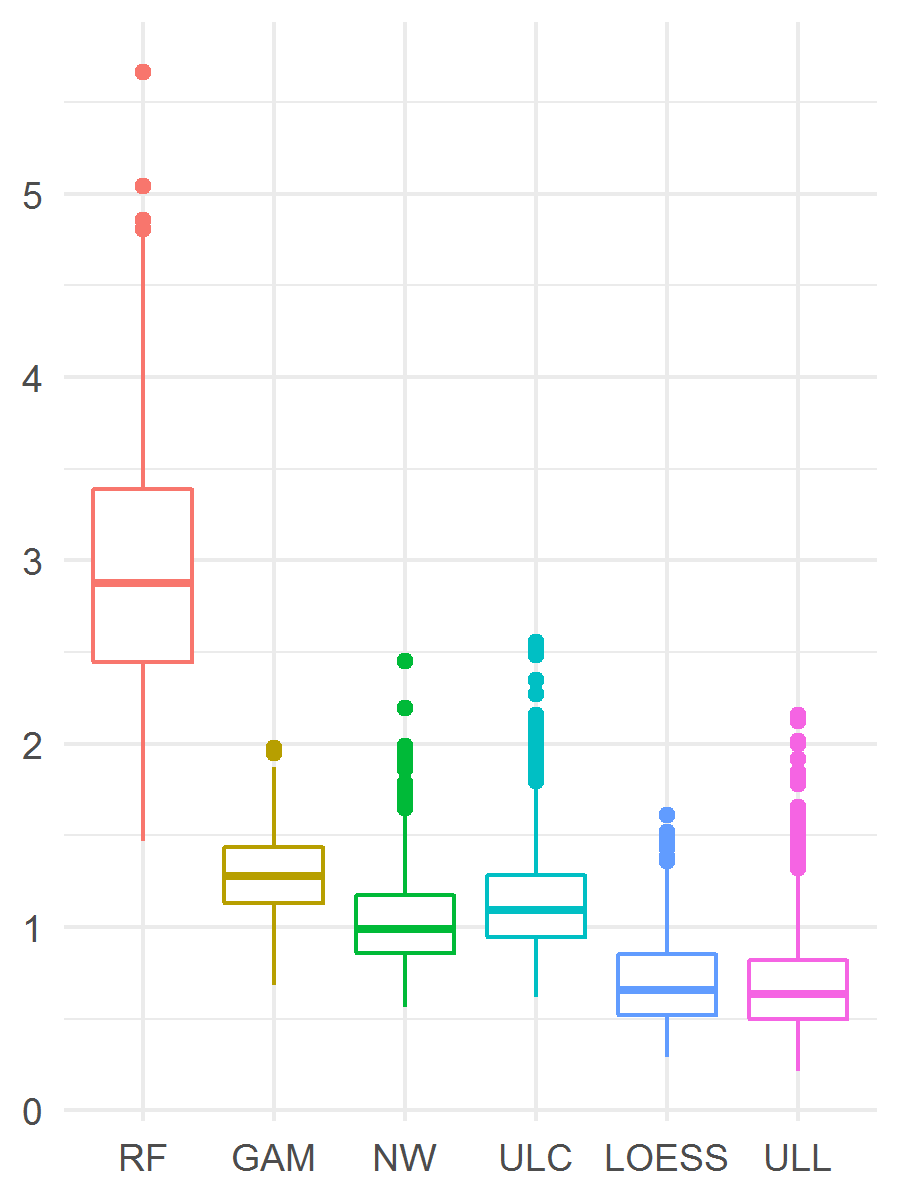} 
  }
  \quad
  \subfigure{
    \includegraphics[width=.45\textwidth]{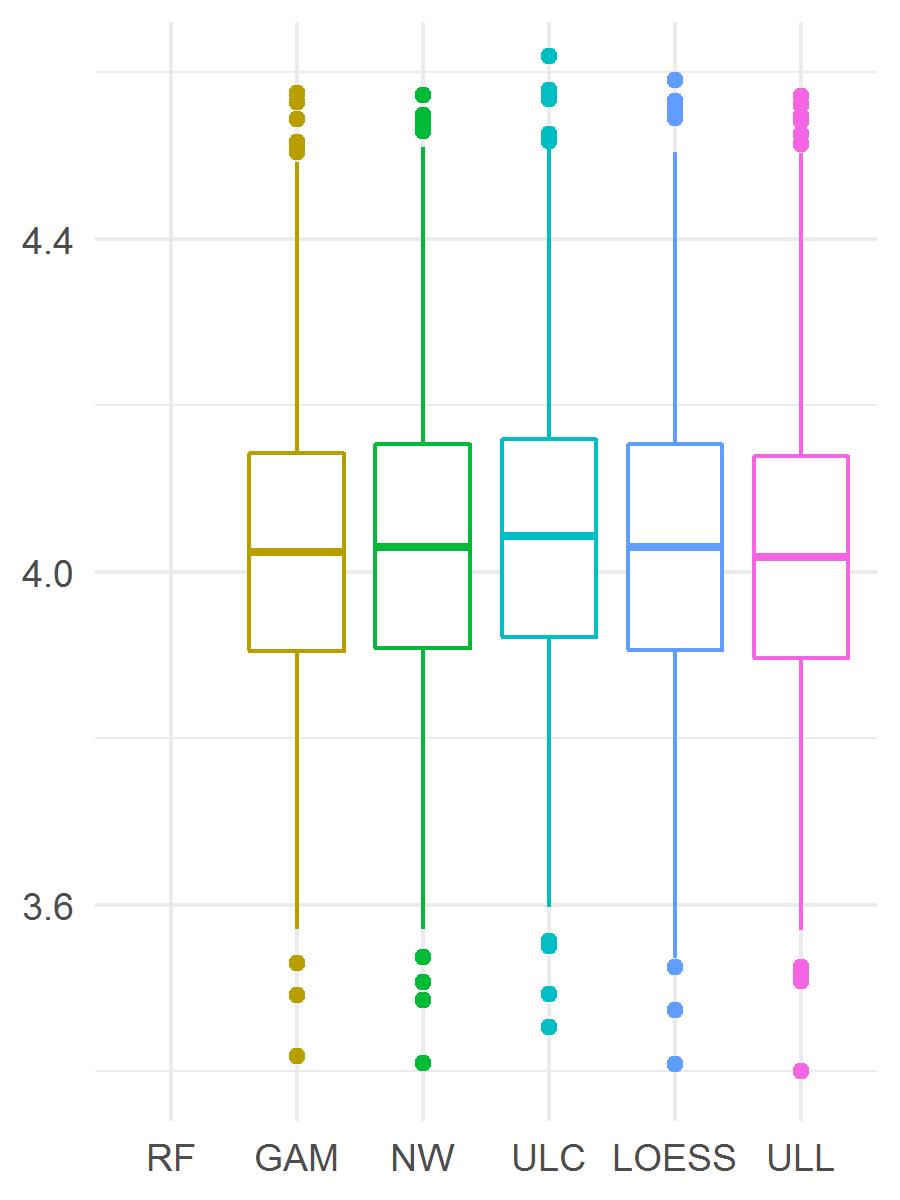} 
  }
  \caption{\footnotesize The maximum (left) and  mean squared (right) errors in Example~2. For the mean squared error, the random forest model performed worse (10.97 (10.55, 11.39)) than the GAM model and the kernel estimators, so the results of the random forest model ``did not fit'' into the plot.}
  \label{fig2}
\end{figure}

The results are presented in Fig.~\ref{fig2}.
For the maximum error, the advantage of the estimators of order 1 (LOESS and (\ref{ull1}))
over the estimators of order 0 (the Nadaraya-Watson and (\ref{ull0})) is noticeable, while the estimator (\ref{ull1})
turns out to be the best of all considered estimators, in particular,
the estimator (\ref{ull1}) performs better than LOESS:
0.6357 (0.4993, 0.8224) vs. 0.6582 (0.5205, 0.8508), $p=0.019$.

For the mean squared error, all models, except random forest and linear regression,
show similar results. Besides, the estimator (\ref{ull1})
turns out to be the best of the considered ones,
although the difference
between estimators (\ref{ull1}) and LOESS is not statistically significant:
4.017 (3.896, 4.139) vs.
4.030 (3.906, 4.154), $p=0.11$.
}
\end{ex}

\bigskip

\begin{ex}
\begin{figure}[!ht]
  \centering
  \includegraphics[width=4in]{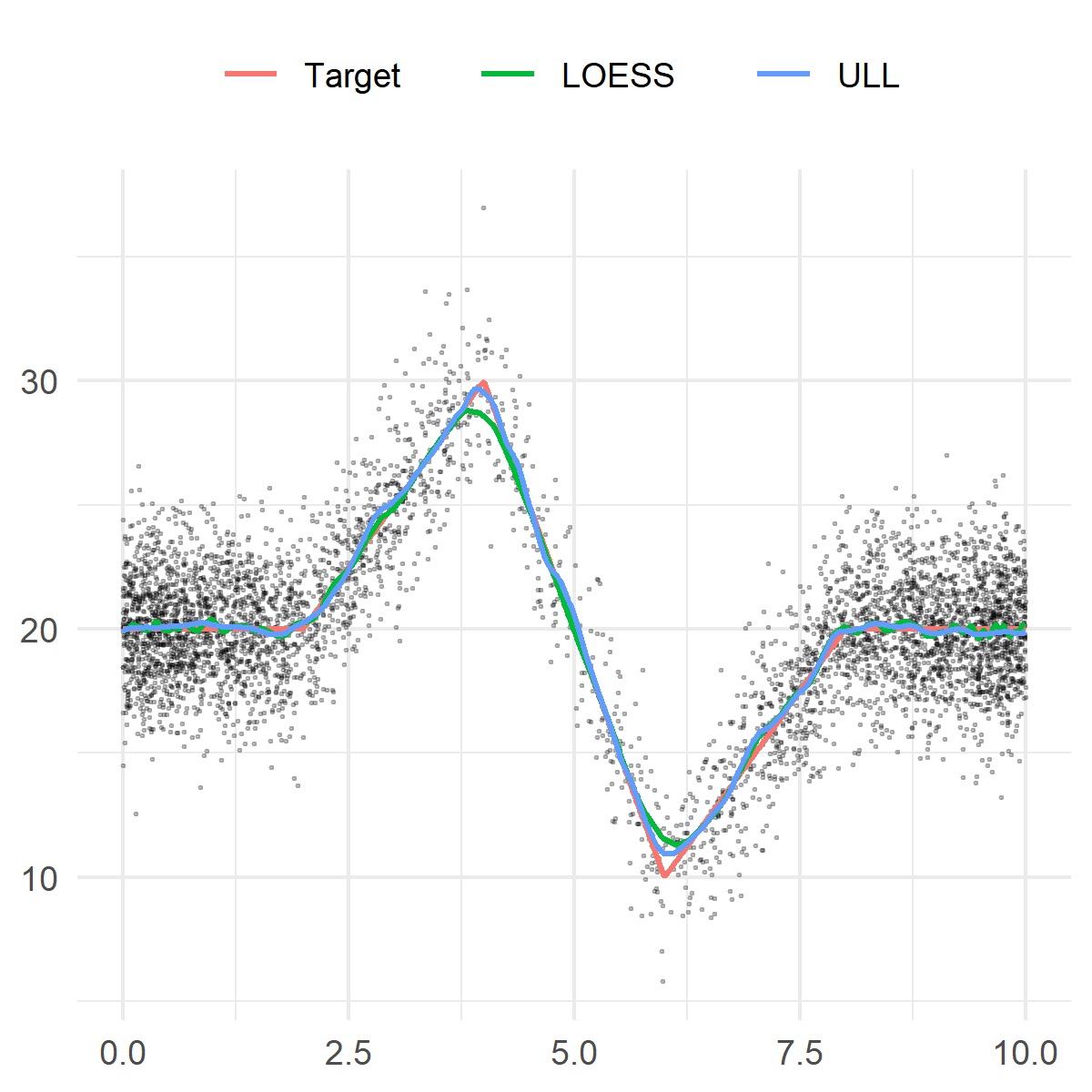}
  \caption{\footnotesize Example 3. Sample observations, target function, and two estimators.}\label{fexam2}
\end{figure}
{\rm
The piecewise linear target function is shown in Fig.~\ref{fexam2}. For the sake of simplicity of presentation, we do not present the formula for the definition of this function.
Here the centered Gaussian noise has the standard deviation $\sigma=2$.
The design points are independent and identically distributed with density
proportional to the function $(z-5)^2+2$, $0\le z \le 10$.

\begin{figure}[!ht]
\centering
  \subfigure{
    \includegraphics[width=.45\textwidth]{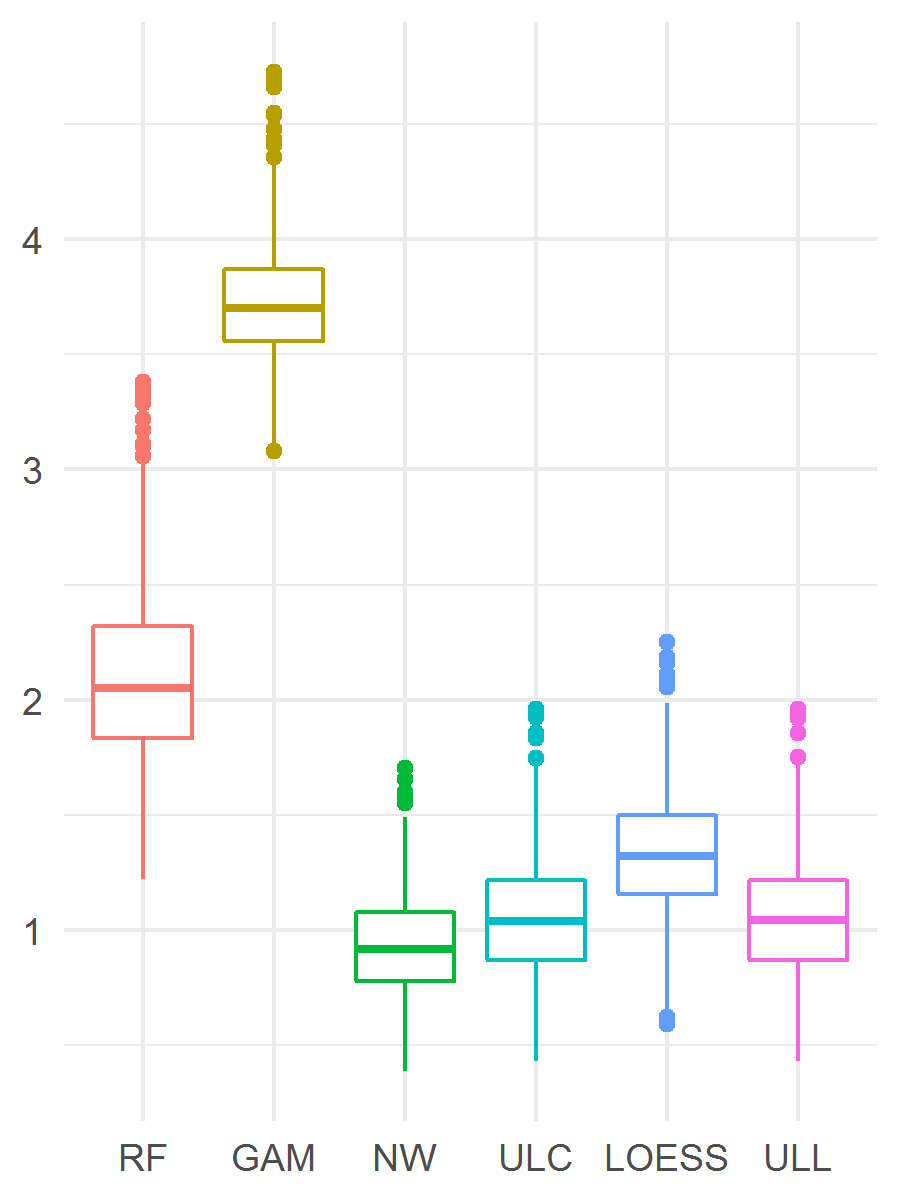} 
  }
  \quad
  \subfigure{
    \includegraphics[width=.45\textwidth]{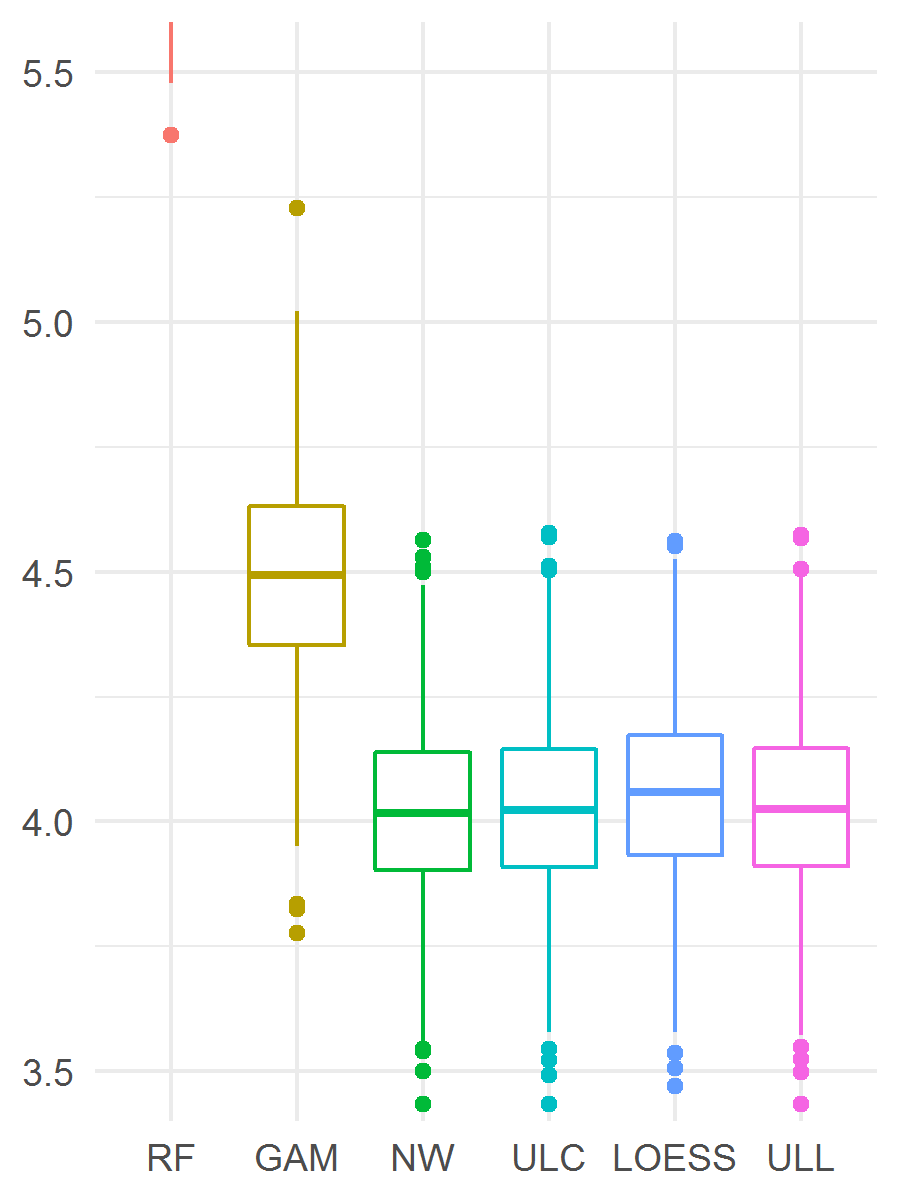} 
  }
  \caption{\footnotesize The maximum (left) and mean squared (right) errors in Example~3. For the mean-squared error, the random forest model
performed worse (6.699 (6.412, 7.046)) than the GAM model and the kernel estimators, so the results of the random forest model ``did not fit'' into the plot.}
\label{rexam2}
\end{figure}

The results are presented in Fig.~\ref{rexam2}.
The Nadaraya-Watson estimator appears to be the best model both for the maximum error and for the mean squared error.
For the both errors,  estimator (\ref{ull1}) is better than  LOESS ($p<0.0001$
for the maximum error, $p=0.0030$ for  the mean squared error).
}
\end{ex}

\bigskip

\begin{ex}
{\rm
In this example, the design points are strongly dependent. We will define them as follows:
$z_i:=s(A_i)$, $i=1,..., n$, where $A$
 is a positive number such that $A/\pi$ is irrational (we chose $A=0.0002$ in this example),
$$s(t):=10 \Big|\sum\nolimits_{k=1}^{100} \eta_k \cos(tk)\Big|
\quad\mbox{ with }\quad \eta_k:={k^{-1}\psi_k}\Big({\sum\nolimits_{j=1}^{100}j^{-1}\psi_j}\Big)^{-1},$$
and $\psi_j$ are independent uniformly distributed on
$[0,1]$ random variables independent of the noise.
It was shown \cite{2020} that the random sequence $s(A_i)$
 is asymptotically everywhere dense on $[0,10]$ with probability~1.

The target function is
$$
f(z) = 0.2\ \big(((z-5)^2+25)\ \cos((z-5)^2/2)\  +\  60\big),
$$
shown in Fig.~\ref{fexam3}.

\begin{figure}[!ht]
  \centering
  \includegraphics[width=4in]{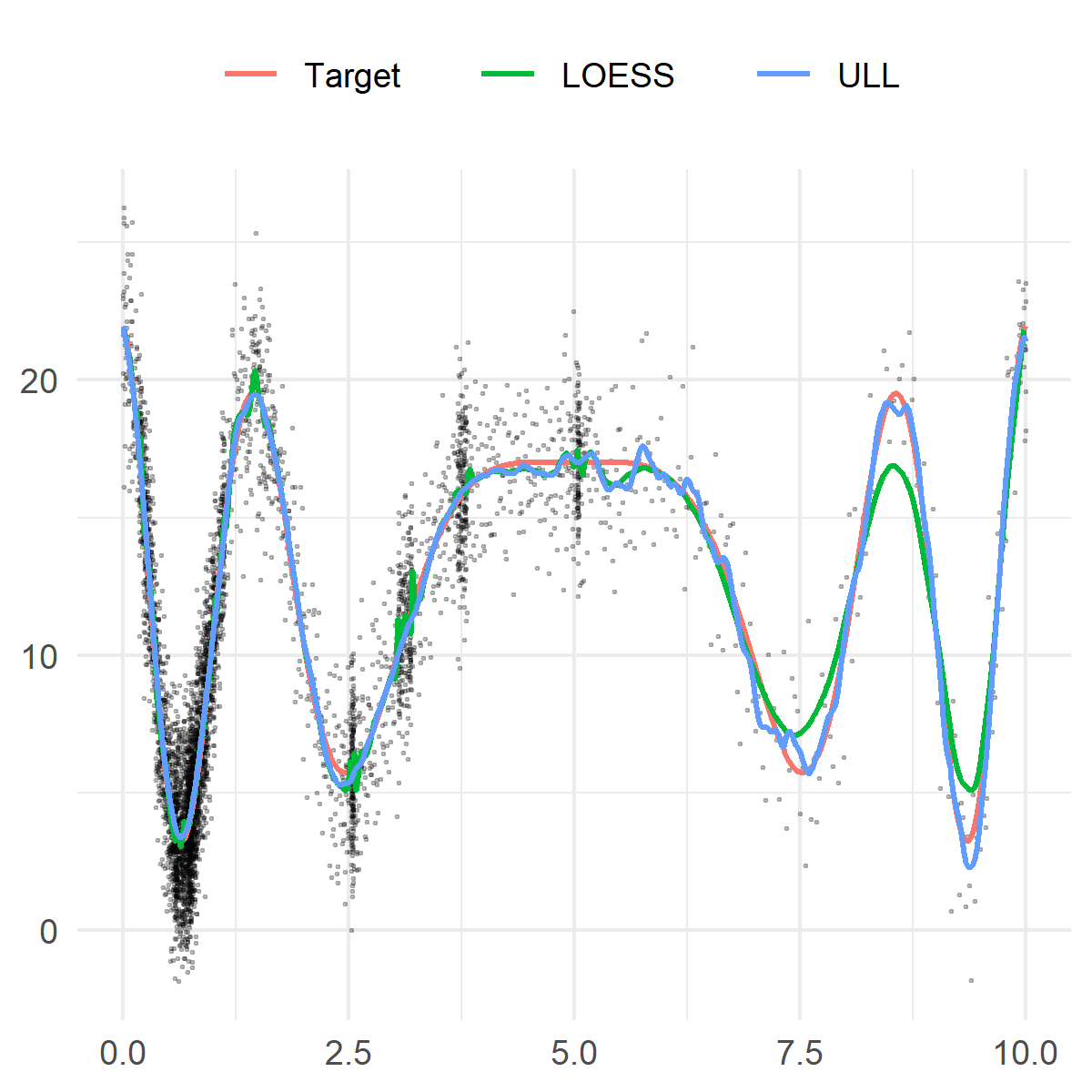}
  \caption{\footnotesize Example 4.  Sample observations, target function, and two estimators.}
\label{fexam3}
\end{figure}

\begin{figure}[!ht]
\centering
  \subfigure{
    \includegraphics[width=.45\textwidth]{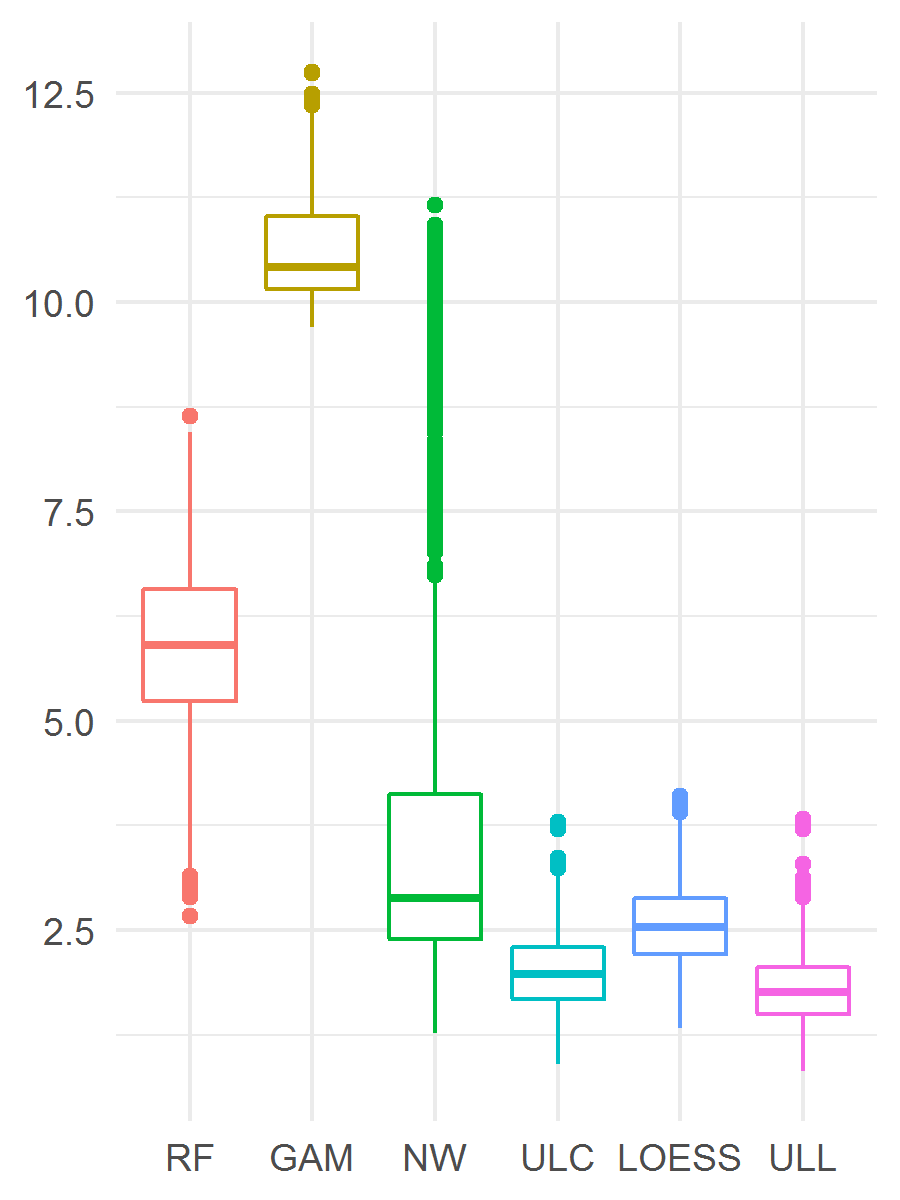} 
  }
  \quad
  \subfigure{
    \includegraphics[width=.45\textwidth]{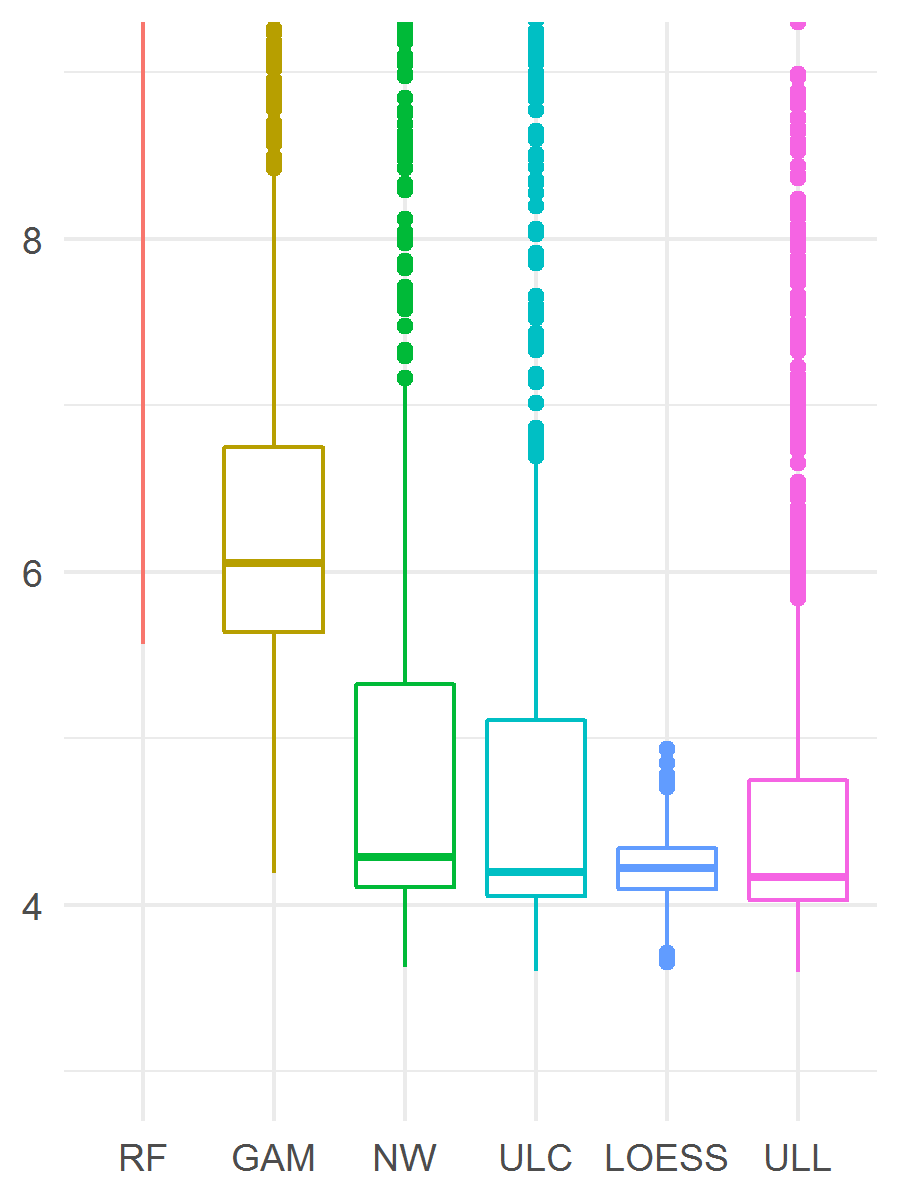} 
  }
  \caption{\footnotesize The maximum (left) and  mean squared (right) errors in  Example~4.
   As before, for the mean squared error, the results of the random forest model
  (13.95 (11.69, 16.18))
 not shown in full on the graph. In addition, the outliers for the GAM, NW, ULC, and ULL estimators are ``cut off'' in this graph.
  } \label{rexam3}
\end{figure}

For maximum error, estimate (\ref{ull1})
turns out to be the best of all the considered estimators. In particular,
estimator (\ref{ull1}) is better than LOESS:
1.757 (1.491, 2.053) vs.
2.538 (2.216, 2.886), $p<0.0001$.

The median mean squared error for
estimator (\ref{ull1}) also turns out to be the smallest of those considered.
In that sense,
estimator (\ref{ull1}) is better than LOESS, but the difference is not significant:
4.166 (4.025, 4.751)
vs. 4.219 (4.096, 4.338), $p=0.92$.
}
\end{ex}

\begin{ex}
{\rm

In this example, the target function was the same as in Example 4.
The difference from the previous example is that
50,000 design points were generated by the same technique, and then 5,000 points of the 50,000 ones were selected.
This allowed us to fill the domain of $f$  with design elements
``more uniformly'' than in the previous example, while preserving the clusters of design points.

\begin{figure}[!ht]
\centering
  \subfigure{
    \includegraphics[width=.45\textwidth]{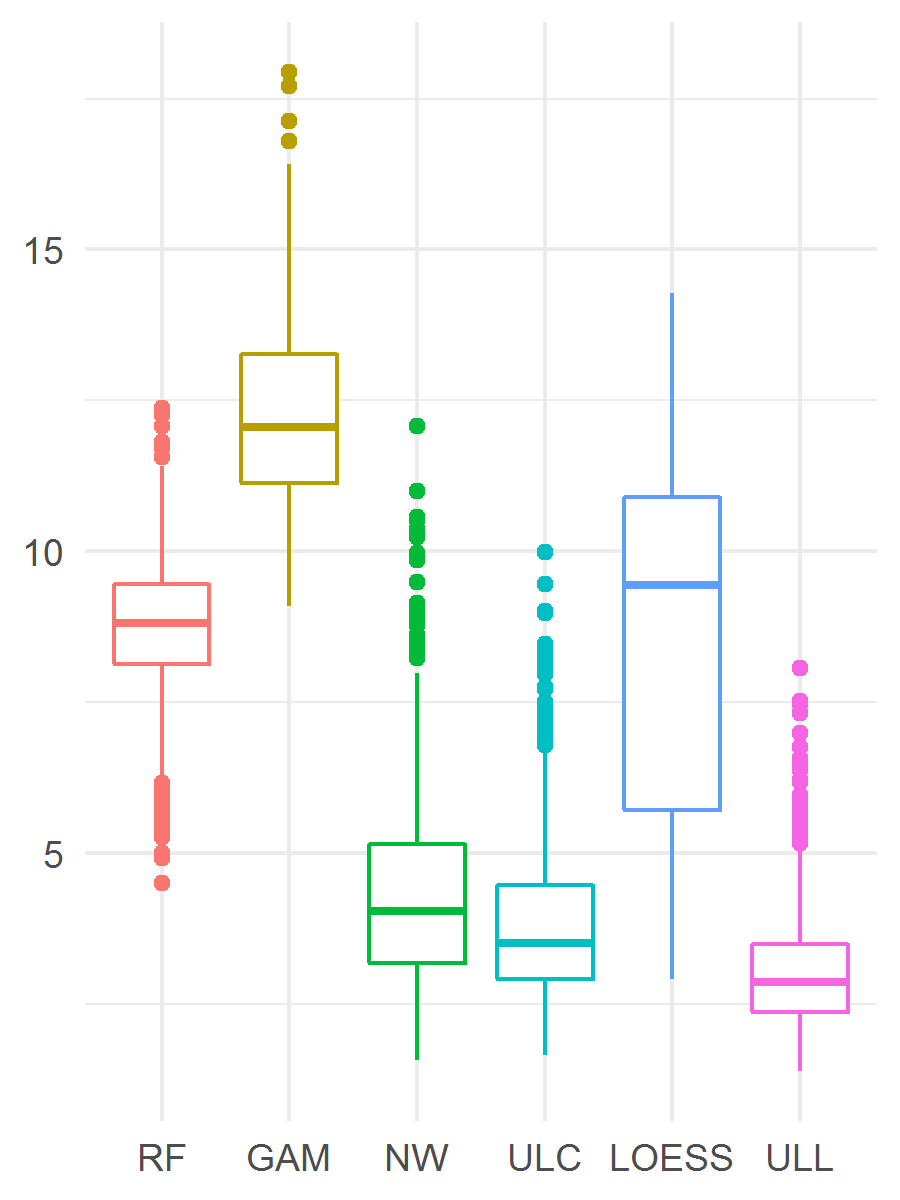} 
  }
  \quad
  \subfigure{
    \includegraphics[width=.45\textwidth]{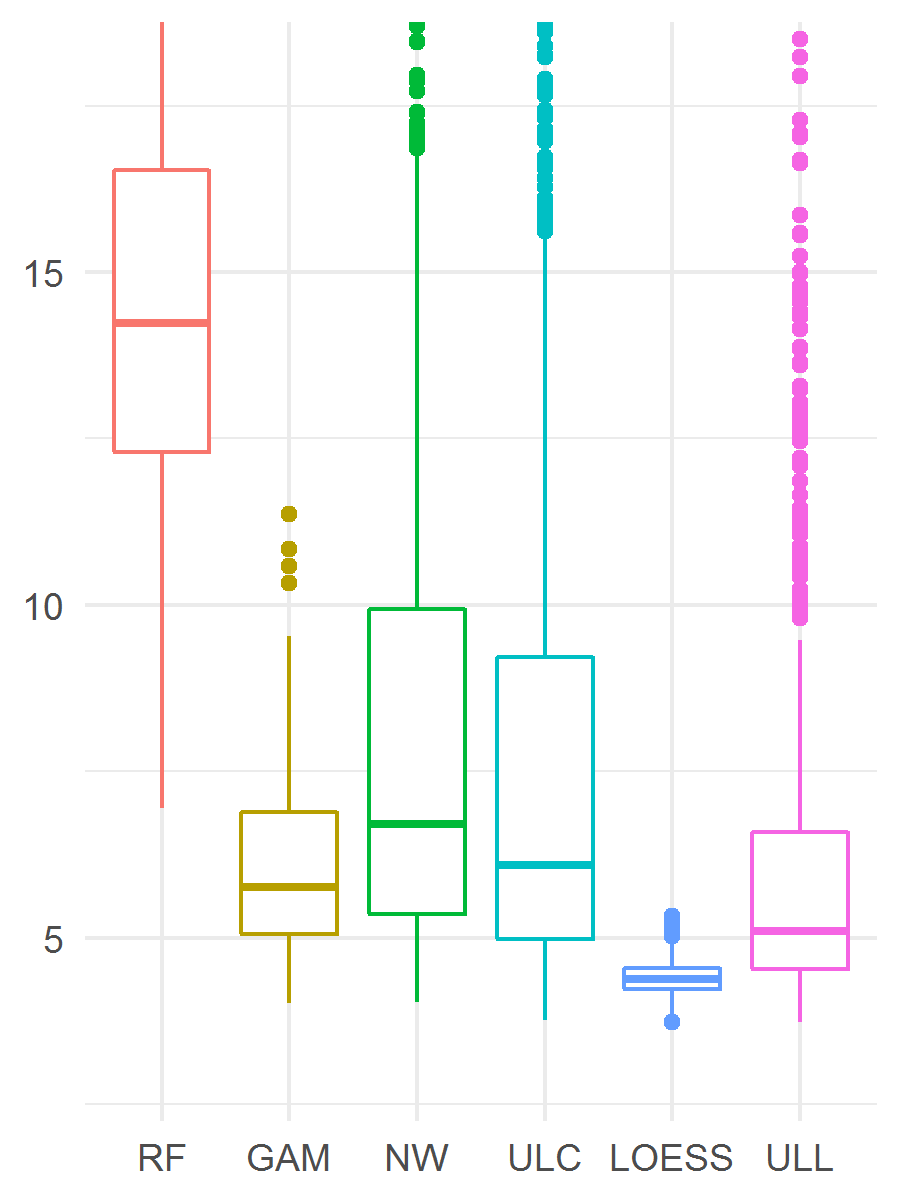} 
  }
  \caption{\footnotesize The maximum (left) and mean-squared (right) errors in Example~5.
     As before, for the mean-squared error, the results of the random forest model
    not shown in full on the graph. In addition, the outliers for the NW, ULC, and ULL  estimators are ``cut off'' in this graph.
  } \label{rexam4}
\end{figure}

For maximum error, estimator (\ref{ull1})
turns out to be the best of all the considered estimators. In particular,
estimator (\ref{ull1}) is better than LOESS:
2.872 (2.369, 3.488) vs.
9.435 (5.719, 10.9), $p<0.0001$.

For the mean squared error, the best estimator is LOESS.
Estimator \eqref{ull1} is worse than LOESS:
5.108 (4.535, 6.597) vs.
4.378 (4.229, 4.541), $p<0.0001$, but it is better than the other estimators considered.
}
\end{ex}

\section{Example of processing real medical data } \label{sec-real}

In this section, we consider an application of the models considered in the previous section
to the data collected in the multicenter study ``Epidemiology of cardiovascular diseases in the regions of the Russian Federation''.
In that study, representative samples  of unorganized male and female populations aged 25--64 years from 13 regions of the Russian Federation were studied.
The study was approved by the Ethics Committees of the three federal centers: State Research Center for Preventive Medicine,
Russian Cardiology Research and Production Complex,  Almazov Federal Medical Research Center. Each participant has
written informed consent for the study. The study was described in detail in \cite{2020-ShD}.

One of the urgent problems of modern medicine is to study the relationship between heart rate (HR) and systolic arterial blood pressure (SBP), especially for low values of the observations.
Therefore we will choose SBP as the outcome, and HR as the predictor.
 The association between these variables was previously estimated to be nonlinear \cite{2020-ShKK}. The general analysis included 6597 participants from 4 regions of the Russian Federation.
The levels of SBP and HR were statistically significantly pairwise different between the selected regions. Thus, the hypothesis of the independence of design points was violated.

In this section, the maximum error cannot be calculated because the exact form of the relationship is unknown, so only the  mean squared error is reported.
The mean squared error was calculated for 1000 random partitions of the entire set of observations into training ($80\%$) and validation ($20\%$) samples.

\begin{figure}[!ht]
\centering
    \includegraphics[width=.45\textwidth]{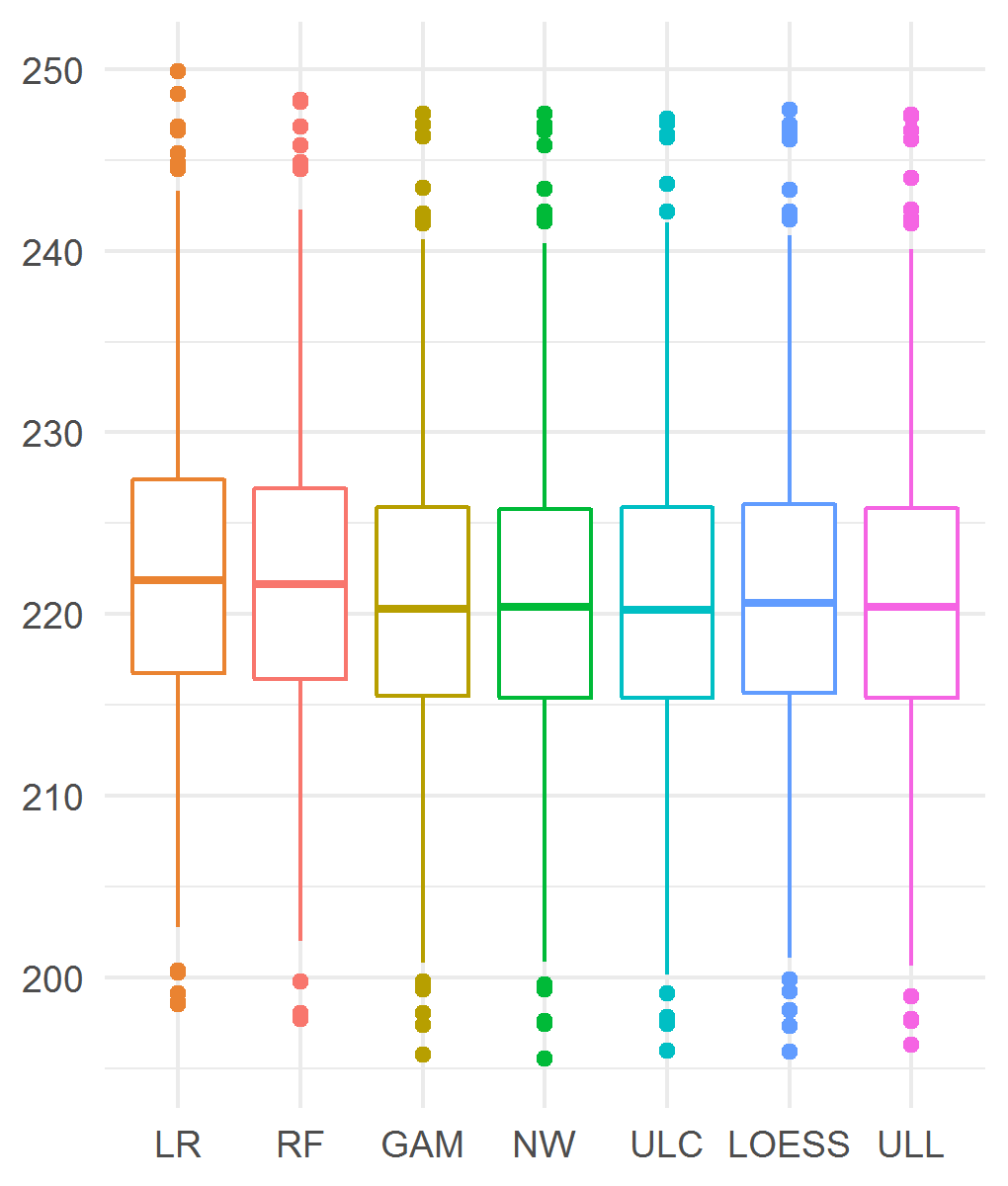} 
  \caption{\footnotesize Mean-squared  prediction error of the dependence of BP from HR.} \label{rexam5}
\end{figure}

the results are presented in Fig.~\ref{rexam5}.
Here the GAM estimator and the kernel estimators showed similar results, which were better than the results
of both the linear regression and random forest.

The best estimator turned out to be \eqref{ull0}, although its difference from
the Nadaraya-Watson estimator
was not statistically significant:
220.2 (215.4, 225.9) vs.
220.4 (215.4, 225.8),
$p=0.91$.
The difference between estimator (\ref{ull1}) and LOESS was not significant too:
220.4 (215.4, 225.9)
vs. 220.6 (215.6, 226.1), $p=0.52$.

\section{Conclusion}

In this paper, for a wide class of nonparametric regression models with a random design,
 universal uniformly consistent kernel estimators are  proposed for an unknown random regression function of a scalar argument. These estimators belong
to the class of local linear estimators.
But in contrast to the vast majority of previously known results, traditional conditions of dependence of design elements
are not needed for the consistency of the new estimators.
The design can be either fixed and not necessarily regular, or random and not necessarily
consisting of independent or weakly dependent random variables. With regard to design elements, the only condition  that is required is the dense filling of the regression function domain with the design points.

Explicit upper bounds are found for the rate of uniform convergence in probability of the new estimators to an unknown random regression function.
The only characteristic explicitly included in these estimators is the maximum spacing statistic of the variational series of design elements,
which requires only the convergence to zero in probability of the maximum spacing as the sample size tends to infinity.
The advantage of this condition over the classical ones
is that it is insensitive to the forms of dependence of the design observations.
Note that this condition is, in fact, necessary, since only when the design densely fills the regression function domain, it is possible to reconstruct the regression function with some accuracy. As a corollary of the main result, we obtain consistent estimators for the mean function of continuous random processes.

In the simulation examples of Section \ref{sec-sim}, the new estimators were compared
with known kernel estimators. In some of the examples,
the new estimators proved to be the most accurate.
In the application to real medical data considered in Section \ref{sec-real},
the accuracy of new estimators was also comparable with that of the best-known kernel estimators.

\section{Proofs}

In this Section, we will prove the assertions stated in Sections 2--4.
 Denote
\begin{equation}\label{beta0}
\beta_{n,i}(t):=\frac{w_{n2}(t)-(t-z_{n:i})w_{n1}(t)}{w_{n0}(t)w_{n2}(t)-w_{n1}^2(t)}.
\end{equation}
Taking into account the relations $X_{ni}=f(z_{n:i})+\varepsilon_{ni}$,\, $i=1,\ldots,n, $ and the identity
\begin{equation}\label{weight}
\sum_{i=1}^n\beta_{n,i}(t)K_{h}(t-z_{n:i})\Delta z_{ni}\equiv 1,
\end{equation}
we obtain the representation
\begin{equation}\label{30}
 \widehat f_{n,h}(t)=f(t)+f(t)I(\delta_n> c_*h)+\widehat r_{n,h}(f,t)+\widehat \nu_{n,h}(t),
\end{equation}
where
\begin{eqnarray*}
\widehat  r_{n,h}(f,t)=I(\delta_n\le c_*h)\sum_{i=1}^n\beta_{n,i}(t)(f(z_{n:i})-f(t))K_{h}(t-z_{n:i})\Delta z_{ni},
\\
\widehat \nu_{n,h}(t)=I(\delta_n\le c_*h)\sum_{i=1}^n\beta_{n,i}(t)K_{h}(t-z_{n:i})\Delta z_{ni}\varepsilon_{ni}.
\end{eqnarray*}
We emphasize that, in view of the properties of the density $ K_h(\cdot) $, the domain of summation in the last two sums as well as in all sums defining the quantities $ w_{nj}(t) $ from (\ref{2}) coincides with the set $ A_{n,h}(t) = \{i: \, | t-z_{n: i}| \le h,\, 1 \le i \le n \} $, which is a crucial point for further analysis.

\begin{lem}\label{lem-2}
For  $h<1/2$, the following equalities are valid$:$
\begin{eqnarray}
\label{17}
\inf\limits_{t\in [0,1]}(w_{0}(t)w_{2}(t)-w_{1}^2(t))= \frac{1}{4}(\kappa_2-\kappa^2_1)h^{2},\quad \inf\limits_{t\in [0,1]}w_{0}(t)=1/2,\\
\label{16-}
\sup\limits_{t\in [0,1]}|w_{j}(t)|=\left(\frac{1}{2}\right)^{j-2[j/2]}\kappa_jh^j, \quad j=0,1,2,3.
\end{eqnarray}
Moreover, on the set of elementary events such that  $\delta_n\leq c_*h$, the following inequalities hold$:$
\begin{eqnarray}
\label{16}
\sup\limits_{t\in [0,1]}|w_{nj}(t)|\le 3Lh^j,\quad
\sup\limits_{t\in [0,1]}|w_{nj}(t)-w_{j}(t)|\le 12L\delta_nh^{j-1},\quad  j=0,1,2,3,\\
\label{18}
\inf\limits_{t\in [0,1]}(w_{n0}(t)w_{n2}(t)-w_{n1}^2(t))\ge \frac{1}{8}(\kappa_2-\kappa^2_1)h^{2},\quad \inf\limits_{t\in [0,1]}w_{n0}(t)\geq1/4,\\
\label{18+1}
\forall t_1, t_2\in [0,1]\quad |w_{nj}(t_2)-w_{nj}(t_1)|\le 18Lh^{j-1}|t_2-t_1|,\quad j=0,1,2.
\end{eqnarray}
\end{lem}
{\it Proof}. Let us prove   (\ref{17}) and (\ref{16-}).
First of all, note that, due to the Cauchy--Bunyakovsky--Schwartz inequality, $ w_{0}(t)w_{2}(t) -w_{1}^2(t) \ge 0 $ for all $ t \in [0,1] $ and this  difference is continuous in $ t $. First, consider the simplest case where $ h \le t \le 1-h $. For such $ t $, after changing the integration variable
in the definition  (\ref{int}) of the quantities $w_{j}(t) $ we have
\begin{equation}\label{identW}
w_{j}(t)=\int\limits_{t-h}^{t+h}(t-z)^jK_{h}(t-z)dz=h^{j}\int\limits_{-1}^{1}v^jK(v)dv,
\end{equation}
i.e., $w_{0}(t)\equiv 1$, $w_{1}(t)\equiv 0$, and $w_{2}(t)\equiv h^{2}\kappa_2$.
In other words,  on the segment $[h,1-h]$, the following identity  is valid$:$
\begin{equation}\label{ident}
w_{0}(t)w_{2}(t)-w_{1}^2(t)\equiv h^{2}\kappa_2.
\end{equation}

We now consider the case $t=\alpha h$ for all $\alpha\in [0,1]$. Then
\begin{equation}\label{zero}
w_{j}(\alpha h)=\int\limits_{0}^{(1+\alpha)h}(\alpha h-z)^jK_{h}(\alpha h-z)dz=h^{j}\kappa_j(\alpha).
\end{equation}
Next, by (\ref{zero}), we obtain
\begin{multline*}
\frac{d}{d\alpha}h^{-2}(w_{0}(\alpha h)w_{2}(\alpha h)-w_{1}^2(\alpha h))=\frac{d}{d\alpha}(\kappa_0(\alpha)\kappa_2(\alpha)-\kappa^2_1(\alpha))\\
=K(\alpha)\left(\alpha^2\int\limits_{-1}^{\alpha}K(v)dv+
\int\limits_{-1}^{\alpha}v^2K(v)dv-2\alpha\int\limits_{-1}^{\alpha}vK(v)dv\right)\ge 0
\end{multline*}
in view of the relation $\int\nolimits_{-1}^{\alpha}vK(v)dv\le 0$  since $K(v)$ is an even function.
Similarly we study the symmetrical case where $t=1-\alpha h$ for all $\alpha\in [0,1]$.
>From here and (\ref{ident}) we obtain the first relation in (\ref{17}):
\begin{equation*}
\inf_{t\in [0,1]}\{w_{0}(t)w_{2}(t)-w_{1}^2(t)\}=w_{0}(0)w_{2}(0)-w^2_{1}(0)=\frac{1}{4}h^{2}(\kappa_2-\kappa^2_1).
\end{equation*}
The second relation in (\ref{17}) directly follows from (\ref{zero}). Moreover, the above-mentioned arguments and the representations
 (\ref{identW}) and  (\ref{zero}) imply (\ref{16-}).

Further, the first estimator in (\ref{16})
is obvious by the above remark about the domain of summation
in the definition of functions $ w_{nj}(t) $,
and the relations
\begin{equation}\label{error-}
\sup\limits_{s\in[0,1]}K(s)\leq L,\quad \sum\limits_{i\in A_{n,h}(t)}\Delta z_{ni}\leq 2h+\delta_n\leq 3h.
\end{equation}
The second estimator in (\ref{16}) immediately follows from the well-known estimate of the error of approximation by Riemann integral sums of the corresponding integrals of smooth functions on a finite closed interval:
\begin{equation}\label{error}
\Big|\sum_{i\in A_{n,h}(t)}g_{t,j}(z_{n:i})\Delta z_{ni}-\int\limits_{z\in [0,1]: |t-z|\leq h}g_{t,j}(z)dz\Big|\le(2h+\delta_n)\delta_n
L_{g_{t,j}},
\end{equation}
where the functions
 $g_{t,j}(z)=(t-z)^jK_{h}(t-z)$,  $j=0,1,2,3$, are defined for all $z\in [0\vee t-h,1\wedge t+h]$, and $L_{g_{t,j}}$ is the Lipschitz constant
 of the function $g_{t,j}(z)$;
It easy to verify that
$\sup_{t\in[0,1]}L_{g_{t,j}}\leq 4Lh^{j-2}$ for all $h\in (0,1/2)$ and $j=0,1,2,3$.
So, on the set of elementary events such that $ \{\delta_n\le c_*h \} $ (recall that $ c_* <1 $), the right-hand side in (\ref{error}) can be replaced
with $12L\delta_nh^{j-1}$.

In addition, taking (\ref{16-}) and (\ref{16})  into account, we obtain
\begin{multline*}
|w_{n0}(t)w_{n2}(t)-w_{0}(t)w_{2}(t)|\\\le w_{n0}(t)|w_{n2}(t)-w_{2}(t)|+w_{2}(t)|w_{n0}(t)-w_{0}(t)|
\le 9L\delta_n(3L+\kappa_2)h,
\end{multline*}
$$|w^2_{n1}(t)-w^2_{1}(t)|\le |w_{n1}(t)-w_{1}(t)|(|w_{n1}(t)|+|w_{1}(t)|)\le 9L\delta_n(3L+\kappa_1/2)h.
$$
Hence follows the estimate
\begin{equation}\label{error2}
|w_{n0}(t)w_{n2}(t)-w_{n1}^2(t)-w_{0}(t)w_{2}(t)+w_{1}^2(t)|\le 9L\delta_n(6L+\kappa_2+\kappa_1/2)h.
\end{equation}
The inequalities in (\ref{18}) follow from (\ref{17}), (\ref{error2}), and the definition of the constant $ c_* $.

To prove (\ref{18+1}), note that
$$w_{nj}(t_2)-w_{nj}(t_1)=\sum\limits_{i=1}^n\left\{(t_2-z_{n:i})^jK_{h}(t_2-z_{n:i})-(t_1-z_{n:i})^jK_{h}(t_1-z_{n:i})\right\}\Delta z_{ni}
$$
$$=\sum\limits_{i\in A_{n,h}(t_1)\cup A_{n,h}(t_2)}\left\{(t_2-z_{n:i})^jK_{h}(t_2-z_{n:i})-(t_1-z_{n:i})^jK_{h}(t_1-z_{n:i})\right\}\Delta z_{ni}
$$
where we can use the estimates
  $|(t_2-z_{n:i})^j-(t_1-z_{n:i})^j|\le 2h^{j-1}|t_2-t_1|$ for $j=0,1,2$, $|t_k-z_{n:i}|\le h$ for $k=1,2$, and also the inequalities
$$|K_{h}(t_2-z_{n:i})-K_{h}(t_1-z_{n:i})|\le Lh^{-2}|t_2-t_1|,$$
\begin{equation}\label{27}
\sum\limits_{i\in A_{n,h}(t_1)\cup A_{n,h}(t_2)}\Delta z_{ni}\le 4h+2\delta_n\le 6h.
\end{equation}

Thus, Lemma \ref{lem-2} is proved. \hfill$\square$

\begin{lem}\label{lem-1}
For any positive $h<1/2$, the following estimate is valid$:$
\begin{eqnarray*}
\sup_{ t\in [0,1]}|\widehat r_{n,h}(f,t)|\le C^*_1\omega_f(h),\quad\mbox{with}\quad C_1^*=C_1\frac{L^2}{\kappa_2-\kappa_1^2}.
\end{eqnarray*}
\end{lem}
{\it Proof}. Without loss of generality, the required estimate can be derived on the set of elementary events determined by the condition $ \delta_n \leq c_*h $. Then the assertion of the lemma follows from the inequality
\begin{multline}\label{rbound}
|\widehat r_{n,h}(f,t)|\le
\frac{\omega_f(h)w_{n2}(t)}{w_{n0}(t)w_{n2}(t)-w_{n1}^2(t)}\sum_{i\in A_{n,h}(t)}K_{h}(t-z_{n:i})\Delta z_{ni}\\
+\frac{\omega_f(h)|w_{n1}(t)|}{w_{n0}(t)w_{n2}(t)-w_{n1}^2(t)}\sum_{i\in A_{n,h}(t)}|t-z_{n:i}|K_{h}(t-z_{n:i})\Delta z_{ni},
\end{multline}
the estimates from (\ref{error-}), and Lemma \ref{lem-2}.
\hfill$\square$

\begin{lem}\label{lem-3}
 For any  $y>0$ and $h<1/2$, on the set of elementary events such that
  $\delta_n\leq c_*h$, the following estimate is valid$:$
\begin{eqnarray*}
{\mathbb P}_{{\cal F}_n}\left (\sup\limits_{t\in [0,1]}|\widehat \nu_{n,h}(t)|>y\right )\le
C_2^*\sigma^2\frac{\delta_n}{h^2y^2},\quad\mbox{with}\quad  C_2^*=C_2\frac{L^4}{(\kappa_2-\kappa_1^2)^2},
\end{eqnarray*}
where the symbol  ${\mathbb P}_{{\cal F}_n}$ denotes the conditionsl probability given the $\sigma$-field ${\cal F}_n$.

\end{lem}
{\it Proof}. Put
\begin{eqnarray}\label{25}
\mu_{n,h}(t) = \sum\limits_{i\in A_{n,h}(t)}h^{-2}\alpha_{n,i}(t)K_{h}\left(t-z_{n:i}\right)\Delta z_{ni}\varepsilon_{ni},
\end{eqnarray}
where $\alpha_{n,i}(t)=w_{n2}(t)-(t-z_{n:i})w_{n1}(t)$,
and note that from Lemma~\ref{lem-2} and the conditions of Lemma \ref{lem-3} it follows that, firstly, $h^{-2}|\alpha_{n,i}(t)|\le 6L$ if only $i\in A_{n,h}(t)$, and secondly,
\begin{eqnarray}\label{19}
\begin{split}
|\widehat \nu_{n,h}(t)|
\le 8(\kappa_2-\kappa_1^2)^{-1}|\mu_{n,h}(t)|.
\end{split}
\end{eqnarray}

The distribution tail of the random variable
$ \sup_{t \in [0,1]} |\mu_{n,h}(t)| $ will be estimated
 by the so-called chaining proposed by A.N. Kolmogorov to estimate the distribution tail of the supremum norm of a stochastic process with almost surely continuous trajectories (see \cite{1956-Ce}). First of all, note that
the set $ [0,1] $ under the supremum sign can be replaced by the set of dyadic rational points
  $${\cal R}=\{j/2^k;\, j=1,\ldots, 2^k-1; k\ge 1\}.$$
  Thus,
\begin{eqnarray*}\sup_{t\in [0,1]}|\mu_{n,h}(t)|=
\sup_{t\in {\cal R} }|\mu_{n,h}(t)|
\le \max_{j = 1,...,2^m - 1}|\mu_{n,h}(j 2^{-m})|\\ +
\sum_{k=m+1}^\infty \max_{j=1,...,2^k-2}
\big|\mu_{n,h}((j+1) 2^{-k})-\mu_{n,h}(j 2^{-k})\big|,
\end{eqnarray*}
where the natural number  $m$ is defined by the equality $m = \lceil|\log_2 h|\rceil$ (here $\lceil a \rceil$ is the minimal natural number greater than or equal to $a$.
One has
\begin{eqnarray}
{\mathbb P}_{{\cal F}_n}\Big(\sup_{t\in [0,1]}|\mu_{n,h}(t)|>y\Big)
\le{\mathbb P}_{{\cal F}_n}\Big(\max_{j = 1,...,2^m - 1}|\mu_{n,h}(j 2^{-m})|
>a_m y\Big)
\notag +\\
 +
\sum_{k=m+1}^\infty {\mathbb P}_{{\cal F}_n}\Big(\max_{j=1,...,2^k-2}
\big|\mu_{n,h}((j+1) 2^{-k})-\mu_{n,h}(j 2^{-k})\big|
>a_k y\Big)
\notag \\
\label{20}
\le\sum_{j =1}^{2^m - 1} {\mathbb P}_{{\cal F}_n}(|\mu_{n,h}(j 2^{-m})|
>a_m y)+
\nonumber\\
+
\sum_{k=m+1}^\infty \,\sum_{j=1}^{2^k-2}{\mathbb P}_{{\cal F}_n}\Big(
\big|\mu_{n,h}((j+1) 2^{-k})-\mu_{n,h}(j 2^{-k})\big|
>a_k y\Big),
\end{eqnarray}
where $a_m,a_{m+1},...$  is a sequence of positive numbers
such that $a_m+a_{m+1}+...=1$.

Let us now estimate each of the terms on the right-hand side of (\ref{20}). Using Markov's inequality for the second moment and the estimates (\ref{error-}), we obtain
\begin{eqnarray}\label{21}
{\mathbb P}_{{\cal F}_n}(|\mu_{n,h}(j 2^{-m})|>a_m y)
\le \frac{(6L)^2}{(a_m y)^{2}}
\sum\limits_{i\in A_{n,h}(j2^{-m})}
 K^2_{h}(j2^{-m} -z_{n:i})(\Delta z_{ni})^2\sigma^2\nonumber\\
 \le  (6L)^2\sigma^2 (a_m y)^{-2} \delta_n (2h+\delta_n)h^{-2}
\le C_3L^2\sigma^2  (a_m y)^{-2} \delta_n h^{-1}.
\end{eqnarray}
Further,
\begin{eqnarray}\label{22}
{\mathbb P}_{{\cal F}_n}\Big(
\big|\mu_{n,h}((j+1) 2^{-k})-\mu_{n,h}(j 2^{-k})\big|
>a_k y\Big)\le
(a_k y)^{-2}h^{-4}\nonumber\\ \times\sum\limits_{i=1}^n {\mathbb E}_{\cal F}\Big(
 \big(\alpha_{n,i}((j+1)2^{-k})K_{h}((j+1) 2^{-k}-z_{n:i})-\alpha_{n,i}(j2^{-k})K_{h}(j 2^{-k}-z_{n:i})\big) \Delta z_{ni}
 \varepsilon_{ni}\Big)^2\nonumber\\
\le\sigma^2(a_k y)^{-2}h^{-4}\nonumber\\
\times\sum\limits_{i=1}^n \Big(
 \alpha_{n,i}((j+1)2^{-k})K_{h}((j+1) 2^{-k}-z_{n:i})-\alpha_{n,i}(j2^{-k})K_{h}(j 2^{-k}-z_{n:i})\Big)^2 (\Delta z_{ni})^2
 \nonumber\\
\le Lh^{-2}|u-v|
\le
C_4\sigma^2 L^4(a_k y)^{-2} 2^{-2k}  \delta_n (4h+2\delta_n)h^{-4}\leq
C_5\sigma^2 L^4(a_k y)^{-2} 2^{-2k}  \delta_nh^{-3}.
\qquad
\end{eqnarray}
Here we took into account that the summation range
in (\ref{22}) coincides with the set
$$\left\{i:\,i\in A_{n,h}((j+1) 2^{-k})\cup A_{n,h}(j2^{-k})\right\},$$
and hence, due to the relation
 $|(j+1)/2^{k}-j/2^k|=2^{-k}\le h$ for  $k>m$, the estimate (\ref{27}) is valid for $t_1=j2^{-k}$ and $t_2=(j+1)2^{-k}$.
Moreover, we used the estimates
$$ \sup_tK_{h}(t)\le Lh^{-1},\quad
| K_{h}(u)-K_{h}(v)|\le Lh^{-2}|u-v|,
$$
and took into account the following inequalities
in the above range of parameter changes (see Lemma \ref{lem-2}):
\begin{eqnarray*}
|\alpha_{n,i}((j+1)2^{-k})-\alpha_{n,i}(j2^{-k})|\le C_{6}Lh2^{-k},\qquad
|\alpha_{n,i}(j2^{-k})|\le C_7Lh^{2},
\\
|\alpha_{n,i}((j+1)2^{-k})K_{h}((j+1) 2^{-k}-z_{n:i})-\alpha_{n,i}(j2^{-k})K_{h}(j 2^{-k}-z_{n:i})|\le C_{8}L2^{-k}.
\end{eqnarray*}
We now obtain from (\ref{20})--(\ref{22}) that
$${\mathbb P}_{{\cal F}_n}\left(\sup_{t\in [0,1]}|\mu_{n,h}(t)|>y\right)
\le C_{9}y^{-2}
\sigma^2 L^4  \delta_nh^{-1}
\left( 2^m a_m^{-2}
+h^{-2} \sum_{k=m+1}^\infty 2^{-k+1} a_k^{-2}
\right).
$$
The optimal sequence $ a_k $ minimizing the right-hand side of this inequality is
 $ a_m = c2^{m/3}$ and $ a_k = ch^{-2/3}
2^{(-k + 1)/3} $ for $ k = m + 1, m + 2, ... $, where $ c $ is defined by the relation
$ a_m + a_{m + 1} + ... = 1 $.
For the indicated sequence, we conclude that
\begin{eqnarray*}
{\mathbb P}_{{\cal F}_n}\left(\sup_{t\in [0,1]}|\mu_{n,h}(t)|>y\right)<
\\
\le C_{10} y^{-2}
\sigma^2L^4   \delta_n h^{-1}
\left( 2^{m/3}
+h^{-2/3} 2^{-m/3} \big(2+2^{1/3}+2^{2/3}\big)
\right)^3\le C_{11}y^{-2}
\sigma^2 L^4  \delta_nh^{-2}.
\end{eqnarray*}
The assertion of the lemma follows from  (\ref{19}).
$\hfill\Box$

{\it Proof of Theorem~$\ref{theor-1}$}. The assertion follows from Lemmas  \ref{lem-1} and \ref{lem-3} if we set
  $$\zeta_n(h)=\sup\limits_{t\in[0,1]}|\widehat \nu_{n,h}(t)|+\sup\limits_{t\in[0,1]}|f(t)|I(\delta_n>c_*h)$$ and take into account the relation
$${\mathbb P}\big( \zeta_n(h)>y,\, \delta_n\le c_*h\big)= {\mathbb E} I\big( \delta_n\le c_*h\big){\mathbb P}_{{\cal F}_n}\big( \zeta_n(h)>y\big),$$ which was required. $\hfill\Box$

To prove Theorem \ref{theor-2} we need the two auxiliary assertions below.

\begin{lem}\label{lem-4}
If the condition {\rm(\ref{51})} is fulfilled then $\lim\nolimits_{\varepsilon\to 0}{\mathbb E}\omega_f(\varepsilon)=0$ and for independent copies
of the a.s. continuous random process   $f(t)$ the following strong law of large numbers is valid$:$
As $N\to \infty$,  then
\begin{equation}\label{55}
\sup\limits_{t\in [0,1]}\big|\overline{f}_{N}(t)-{\mathbb E}f(t)\big|\stackrel{p}\to 0, \quad \mbox{where}\quad \overline f_{N}(t)=N^{-1}\sum\limits_{j=1}^Nf_j(t).
\end{equation}
\end{lem}
{\it Proof}. The first assertion of the lemma follows from (\ref{51}) and Lebesgue's dominated convergence theorem. We put
$$\omega_{\overline f_N}(\varepsilon)=\sup\limits_{t,s:|t-s|\leq\varepsilon}\big|\overline{f}_{N}(t)-\overline{f}_{N}(s)\big|,  \quad
\omega_{\mathbb{E} f}(\varepsilon)=\sup\limits_{t,s:|t-s|\leq\varepsilon}\big|\mathbb{E}{f}(t)-\mathbb{E}{f}(s)\big|.$$
For any fixed $k>0$ and $i=0,\ldots,k$, one has
\begin{eqnarray} \label{56}
\sup\limits_{t\in [0,1]}\big|\overline{f}_{N}(t)-{\mathbb E}f(t)\big|\leq \max\limits_{0\leq i\leq k}\left|\overline{f}_{N}\big(i/k\big)-{\mathbb E}f\big(i/k\big)\right|+\nonumber\\+
\max\limits_{1\leq i\leq k}\sup\limits_{(i-1)/k\leq t\leq i/k}\left|\overline{f}_{N}(t)-\overline{f}_{N}\big(i/k\big)\right|+
\max\limits_{1\leq i\leq k}\sup\limits_{(i-1)/k\leq t\leq i/k}\left|\mathbb{E}{f}(t)-\mathbb{E}{f}\big(i/k\big)\right|\leq\nonumber\\
\leq \max\limits_{0\leq i\leq k}\left|\overline{f}_{N}\big(i/k\big)-{\mathbb E}f\big(i/k\big)\right|+\omega_{\overline f_N}\left({1}/{k}\right)+\omega_{\mathbb{E} f}\left({1}/{k}\right).
\end{eqnarray}
Put  $\omega_{ f_j}(\varepsilon)=\sup\limits_{t,s:|t-s|\leq\varepsilon}\big|{f}_j(t)-{f}_j(s)\big|$ and note that
$\omega_{\mathbb{E}f}(\varepsilon)\leq {\mathbb E}\omega_f(\varepsilon)$, and as $N\to\infty$,
$$\overline{f}_{N}\big(i/k\big)\stackrel{p}{\to}{\mathbb E}f\big(i/k\big), \quad \omega_{\overline f_N}(\varepsilon)\leq\frac{1}{N}\sum\limits_{j=1}^N\omega_{f_j}(\varepsilon)\stackrel{p}{\to}{\mathbb E}\omega_f(\varepsilon).$$
Therefore, the right-hand side in (\ref{56}) does not exceed $ {\mathbb E} \omega_f\left(1/k\right) + o_p(1) $ and
by the arbitrariness of $ k $ and the first statement of the lemma, the relation (\ref{55}) is proved.
\hfill$\square$

\begin{lem}\label{lem-5}
Under the conditions of Theorem~{\rm \ref{theor-2}} the following limit relation holds$:$
\begin{equation}\label{57}
\frac{1}{N}\sum_{j=1}^N\Delta_{n,h,j}\stackrel{p}\to 0,
 \quad \mbox{where}\quad \Delta_{n,h,j}=\sup_{t\in [0,1]}| f^*_{n,h,j}(t)-f_j(t)|.
\end{equation}
\end{lem}
{\it Proof}.
Let  the sequences  $h=h_n\to 0$ and  $N=N_n\to\infty$ be such that condition  (\ref{52}).
Introduce the event
$B_{n,h,j}=\{\delta_{n,j}\le c_*h\}$, where $j=1,\ldots,N$.
For any positive $\nu$ one has
\begin{equation}\label{58}
{\mathbb P}\left \{\frac{1}{N}\sum_{j=1}^N\Delta_{n,h,j}>\nu\right\}\le {\mathbb P}\left \{\frac{1}{N}
\sum_{j=1}^N\Delta_{n,h,j}I(B_{n,h,j})>\nu\right\}+N{\mathbb P}(\overline{B_{n,h,1}}).
\end{equation}
Next, from Theorem~\ref{theor-1} we obtain
\begin{eqnarray*}{\mathbb E}\Delta_{n,h,j}I(B_{n,h,j})\le C_1^*{\mathbb E}\omega_f(h)+
\int\limits_0^{\infty}{\mathbb P}\left( \zeta_n(h)>y,\, \delta_n\le c_*h\right)dy\le\\
\le C_1^*{\mathbb E}\omega_f(h)+h^{-1}({\mathbb E}\delta_n)^{1/2}+
\int\limits_{h^{-1}({\mathbb E}\delta_n)^{1/2}}^{\infty}{\mathbb P}\left( \zeta_n(h)>y,\, \delta_n\le  c_*h\right)dy
\le \\
\le C_1^*{\mathbb E}\omega_f(h)+(1+C_2^*\sigma^2)h^{-1}({\mathbb E}\delta_n)^{1/2}.
\end{eqnarray*}
 To complete the proof of the lemma, it remains for the first probability on the right-hand side of (\ref{58}) to apply Markov's inequality, use the last estimate, limit relations (\ref{52}), and the first statement of Lemma \ref{lem-4}.
\hfill$\square$

{\it The proof} of Theorem \ref{theor-2} follows from Lemmas \ref{lem-4} and \ref{lem-5}.
\hfill$\square$

{\it Proof of Proposition }\ref{predl-4}. For the estimator $f_{n,h}^*(t)$ defined in (\ref{est4}), we need the following representation:
 \begin{equation}\label{30+}
  f_{n,h}^*(t)=f(t)+ r_{n,h}^*(f,t)+\nu_{n,h}^*(t),
\end{equation}
where
\begin{eqnarray*}
 r_{n,h}^*(f,t)=w_{n0}^{-1}(t)\sum_{i=1}^n(f(z_{n:i})-f(t))K_{h}(t-z_{n:i})\Delta z_{ni},
\\
\nu_{n,h}^*(t)=w_{n0}^{-1}(t)\sum_{i=1}^nK_{h}(t-z_{n:i})\Delta z_{ni}\varepsilon_{ni}.
\end{eqnarray*}

In view of the representations (\ref{30}) and (\ref{30+}), we obtain
\begin{eqnarray}\label{newbias}
{\rm Bias} \widehat f_{n,h}(t)={\mathbb E}\widehat r_{n,h}(f,t)+f(t){\mathbb P}(\delta_n>c_*h)\nonumber\\
=\sum_{i=1}^n{\mathbb E}\{I(\delta_n\le c_*h)\beta_{n,i}(t)(f(z_{n:i})-f(t))K_{h}(t-z_{n:i})\Delta z_{ni}\}+f(t){\mathbb P}(\delta_n>c_*h),
\\
\label{newbias+}
{\rm Bias} f_{n,h}^*(t)={\mathbb E} r_{n,h}^*(f,t)
\nonumber\\ =\sum_{i=1}^n{\mathbb E}\{I(\delta_n\le c_*h)w_{n0}^{-1}(t)(f(z_{n:i})-f(t))K_{h}(t-z_{n:i})\Delta z_{ni}\}+\tau_n,
\end{eqnarray}
where $|\tau_n|\leq\omega_f(h){\mathbb P}(\delta_n>c_*h)$.
Further, it follows from Lemma \ref{lem-2} that, under the condition $\delta_n\le c_*h$, for any point $t\in [h, 1-h]$ one has
\begin{equation}\label{beta}
\sup\limits_{i\in A_{n,h}(t)}|\beta_{n,i}(t)-w_{n0}^{-1}(t)|\le C_5^*\delta_nh^{-1}.
\end{equation}
When deriving the relation (\ref{beta}), we also took into account that
 $w_{0}(t)=1$ and $w_{1}(t)= 0$ for all $t\in [h,1-h]$ (see the proof of Lemma \ref{lem-2}).
Now, reckoning with the relations
 (\ref{error-}), (\ref{newbias}), (\ref{newbias+}),  (\ref{beta}), and Lemma \ref{lem-2}, it easy to derive the first assertion of the lemma since
\begin{multline}\label{difbias}
|{\rm Bias} \widehat f_{n,h}(t)-{\rm Bias} f^*_{n,h}(t)|
\le C_5^*h^{-1}\omega_f(h){\mathbb E}\left\{\delta_nI(\delta_n\le c_*h)\sum_{i=1}^nK_{h}(t-z_{n:i})\Delta z_{ni}\right\}\\
+(|f(t)|+\omega_f(h)){\mathbb P}(\delta_n>c_*h)\le C_6^*\omega_f(h)h^{-1}{{\mathbb E}\delta_n}+(|f(t)|+\omega_f(h)){\mathbb P}(\delta_n>c_*h).
\end{multline}

To prove the second assertion, first of all, note that
\begin{multline*}
{\mathbb Var}\widehat f_{n,h}(t)={\mathbb Var}\widehat \nu_{n,h}(t)+{\mathbb Var}\left(\widehat r_{n,h}(f,t)+f(t)I(\delta_n>c_*h)\right)\\
={\mathbb Var}\widehat \nu_{n,h}(t)+{\mathbb Var}\widehat r_{n,h}(f,t)+f^2(t){\mathbb P}(\delta_n>c_*h){\mathbb P}(\delta_n\le c_*h),
\end{multline*}
$${\mathbb Var}f^*_{n,h}(t)={\mathbb Var}\nu^*_{n,h}(t)+{\mathbb Var}r^*_{n,h}(f,t).
$$
Thus, we need to compare the two variances on the right-hand side of the first equality with the corresponding variances of the second one.
Using (\ref{error-}) and (\ref{beta}), we get
\begin{multline*}
|{\mathbb Var}\widehat \nu_{n,h}(t)-{\mathbb Var}\nu^*_{n,h}(t)|\le\sigma^2\left|{\mathbb E}\sum_{i=1}^nI(\delta_n\le c_*h)(\beta^2_{n,i}(t)-w_{n0}^{-2}(t))K^2_{h}(t-z_{n:i})(\Delta z_{ni})^2\right|\\
+\sigma^2{\mathbb P}(\delta_n>c_*h)
\le C_7^*\sigma^2h^{-1}{\mathbb E}\left\{\delta_nI(\delta_n\le c_*h)\sum_{i=1}^nhK^2_{h}(t-z_{n:i})\Delta z_{ni}\right\}\\
+ \sigma^2{\mathbb P}(\delta_n>c_*h)
\le C_8^*\sigma^2h^{-1}{\mathbb E}\delta_n;
\end{multline*}
when deriving this estimate, we took into account that
$$\sum_{i=1}^nw_{n0}^{-2}(t))K^2_{h}(t-z_{n:i})(\Delta z_{ni})^2\le 1.
$$

To estimate the difference $|{\mathbb Var}\widehat r_{n,h}(t)-{\mathbb Var}r^*_{n,h}(t)|$, note that the bound
$C_9^*\overline f^2h^{-1}{\mathbb E}\delta_n$
for the modulus of the difference between the squares of the displacements of the random variables $ \widehat r_{n,h}(f,t) $ and $ r^*_{n,h}(f,t) $ is essentially
contained in (\ref{rbound}) and (\ref{difbias}).
 Estimation of the difference of the second moments of the specified random variables
 is done similarly with (\ref{error-}), (\ref{beta}), and (\ref{difbias}):
\begin{eqnarray*}|{\mathbb E}\widehat r^2_{n,h}(f,t)- {\mathbb E}r^{*2}_{n,h}(f,t)|\le {\mathbb E}|\widehat r_{n,h}(f,t)-r^{*}_{n,h}(f,t)|
|\widehat r_{n,h}(f,t)+r^{*}_{n,h}(f,t)|\le C_{10}^*\overline f^2h^{-1}{{\mathbb E}\delta_n},
\end{eqnarray*}
which completes the proof. \hfill$\square$

{\it Proof of Proposition} \ref{predl-5}.
>From the definition of $\beta_{n,i}(t)$ in (\ref{beta0}) it follows that, for any $t\in [0,1]$,
\begin{eqnarray*}\sum_{i=1}^n\beta_{n,i}(t)(z_{n:i}-t)K_{h}(t-z_{n:i})\Delta z_{ni}=0,\\
\sum_{i=1}^n\beta_{n,i}(t)(z_{n:i}-t)^2K_{h}(t-z_{n:i})\Delta z_{ni}=D^{-1}_n(t)(w^2_{n2}(t)-w_{n3}(t)w_{n1}(t))=:B_n(t),
\end{eqnarray*}
where $D_n(t):=w_{n0}(t)w_{n2}(t)-w^2_{n1}(t)$.
Expanding the function $ f(\cdot) $ by the Taylor formula in a neighborhood of the point $ t $ (up to the second derivative), from the above identities we obtain, using (\ref{beta0}),  (\ref{newbias}),
and lemma \ref{lem-2}, that for any point $ t $ we have
\begin{multline}\label{newbias1}
{\rm Bias} \widehat f_{n,h}(t)
={\mathbb E}I(\delta_n\le c_*h)\sum_{i=1}^n\{\beta_{n,i}(t)(f(z_{n:i})-f(t))K_{h}(t-z_{n:i})\Delta z_{ni}\}+f(t){\mathbb P}(\delta_n>c_*h)\\
=\frac{f''(t)}{2}{\mathbb E}I(\delta_n\le c_*h)B_{n}(t)+f(t){\mathbb P}(\delta_n>c_*h)+o(h^2)\\
=\frac{f''(t)}{2}B_{0}(t)+O({\mathbb E}\delta_n/h)+o(h^2);
\end{multline}
moreover, the $ O $- and $ o $-symbols on the right-hand side of (\ref{newbias1}) are uniform in $ t $. Note that $ B_0(t) = O(h^2) $ holds for any $t$.

Next, since for  $j=1,2$ we have $|w_{j}(t)|w_{0}^{-1}(t)\le h^j$ and  $|w_{nj}(t)|w_{n0}^{-1}(t)\le h^j$ for all natural $n$, the following asymptotic representation holds:
\begin{multline}
{\rm Bias} f^*_{n,h}(t)
=\sum_{i=1}^n{\mathbb E}w_{n0}^{-1}(t)(f(z_{n:i})-f(t))K_{h}(t-z_{n:i})\Delta z_{ni}\\
=-f'(t){\mathbb E}\frac{w_{n1}(t)}{w_{n0}(t)}I(\delta_n\le c_* h)+\frac{f''(t)}{2}{\mathbb E}\frac{w_{n2}(t)}{w_{n0}(t)}
I(\delta_n\le c_* h)+ O(h{\mathbb  P}(\delta_n> c_* h))+o(h^2)\\
=-f'(t)\frac{w_{1}(t)}{w_{0}(t)}+\frac{f''(t)}{2}\frac{w_{2}(t)}{w_{0}(t)}
+ O({\mathbb E}\delta_n)+o(h^2).
\end{multline}

 {\it Proof of Corollary $\ref{equiv}$}.
 Without loss of generality, we can assume that $ t \in [h,1-h] $. Then, as noted in the proof of Lemma \ref{lem-2}, for the indicated $ t $, one has
 $w_0(t)=1$, $w_1(t)=0$, and $w_2(t)=\kappa_2h^2$, i.e., $B_0(t)=\kappa_2h^2$. \hfill$\square$

{\it Proof of Corollary  $\ref{atzero}$}. This assertion follows from Proposition
\ref{predl-5} and (\ref{zero}). \hfill$\square$

\section*{Acknowledgments} Yu. Linke, I. Borisov, and P. Ruzankin were supported within the framework of the state contract of the Sobolev Institute of Mathematics, project FWNF-2022-0009.


\end{document}